# Deterministic Integral and Ordinary Differential Equations over Irregular Paths

Yevgeniy Guseynov

Independent research
gyevg@yahoo.com

**Abstract**

We define a deterministic integral with respect to irregular paths as a limit of standard line integrals and completely describe a class of all paths for which this integral exists for functions with Hölder exponent in the range of (0,1]. With the developed integral calculus, we study the existence, uniqueness and continuity of solution of time- and path-dependent ordinary differential equations driven by irregular paths in traditional Hölder spaces. These results can be viewed as a supplement to the Young-Lyons-Gubinelli integration theory, in particular, for Hölder exponents less than 1/3.

## 1. Introduction

We define a deterministic integral with respect to irregular paths not necessarily of bounded variation for $\beta$-Hölder functions with $0 < \beta \le 1$. In this topic Young (1936) extended the classes of functions for which the Reimann-Stieltjes integral exists to functions of bounded p-variation. For Hölder functions, a classical result of Young states that the Reimann-Stieltjes integral is well defined if the integrand $f$ is $\beta$-Hölder continuous and the integrator $g$ is $\alpha$-Hölder continuous with $\alpha + \beta > 1$. The $\alpha$-Hölder functions, $0 < \alpha < 1$, provide a measure of paths irregularity since there is a wide class of such functions that are not of bounded variation (see, e.g., Wiener (1930), Paley et al. (1933), Shidfar and Sabetfakhri (1986)).
Despite several later definitions of deterministic integrals over irregular curves for continuous functions (see Guseynov (2016) and references therein), Young's result was not broadened.
The restriction $\alpha + \beta > 1$ does not allow one to apply the Young integral in many typical cases, for example, the existence of the integral

$$\int_0^1 x(t)\, dx(t),$$

when $x$ is a Brownian motion sample path (BMSP). A suitable stochastic integral was defined by Itô (1944) on a probability space $(\Omega, \Sigma, P)$. This stochastic integral and Itô's (1951) formula play a fundamental role in stochastic calculus. Moreover, in many applications, stochastic processes are observed as their realizations or sample paths, but the Itô integral is defined nonconstructively as a limit in $L_2(\Omega \times R_+)$, which does not allow one to evaluate the integral over BMSPs. There are also several approaches to pathwise constructions of the Itô stochastic integral that provide the methods to approximate it in some topology related to a probability space by the sequences of the convenient stochastic processes, see, e.g., Wong and Zakai (1965), Willinger and Taqque (1989) and references therein. All these approaches are using the probability theory to prove the convergence of pathwise stochastic integrals to the Itô or Stratonovich integrals.

Despite such success in the stochastic integration there is still the need in the deterministic integral over irregular paths beyond the Young integration with the concept of uniform estimates of the processes, in particular, for nonlinear dynamical systems. Such approach was developed in the pioneering work of Lyons (1998), (2002) for $\alpha$ −Hölder paths, $\alpha > 1/3$, with his rough paths integration theory based on a given values of higher indefinite iterated integrals which can be constructed as a solution of an algebraic problem or using some stochastic integration procedure. This theory was extended by Gubinelli (2004) and many others, see recent book Friz and Hairer (2020) and references therein.
In the scale of Hölder classes, for almost all with respect to the Wiener measure (W-a.a.) BMSPs the Young integral does not exist for all functions from a $\beta$-Hölder class if $\beta \le \frac{1}{2}$, see Corollary 2.2. This led Gubinelli (2004) to introduce a subclass of $\beta$-Hölder functions whose increments are controlled by a given rough path of $\alpha$-Hölder function, $\alpha \geq \beta > 1/3$, where an integral can be defined and the problem of the existence, uniqueness and continuity of solution of differential equations driven by such paths can be studied. Thus, one of the main problems solved by Lyons-Gubinelli theory is to describe for each class of rough controlling paths a wide class of functions (Gubinelli class) for which an integration process can be defined for each such controlling path.

In the present paper we consider an opposite problem: for the given class of $\beta$-Hölder functions describe a subclass of $\alpha$-Hölder paths, $1 \geq \alpha, \beta > 0$, where the introduced integral with respect to irregular paths exists for all functions from given $\beta$-Hölder class. As we show later, for any $\beta > 0$ such a class of irregular path includes many important in dynamical systems highly oscillatory inputs, see Example 4.1.

A multidimensional integral that we use in Section 3 to study ordinary differential equations driven by irregular paths (ODEIP) is defined by the one-dimensional integral

$$\int_0^t f(s, g(s)) \, dg(s),$$

where $t \in R_+$ and restrictions on $f$ and $g$ are described in Definition 2.2. It could be presented as a line integral

$$\int_{\gamma(\widetilde{0}),\gamma(t)} \omega,$$

where the contour $\gamma$ is given by a parametric presentation $\gamma(t) = (t, g(t)), t \in [0,1]$ in the plane $(t, x)$ and the differential form $\omega(t, x) = f(t, x) dx$. With deterministic integrands and integrators, a necessary and sufficient condition for its existence was proven by Guseynov (2016) in terms of the metric characteristics of a region boundary in $R^d, d \geq 2$, particularly for nonrectifiable Jordan curves in the plane. In the present paper, given the geometry of paths $\gamma$, we adjusted the (2016) integral to introduce the new definition which making it possible presenting the integral through a special approximation sequence of segments in the plane that converges to $\gamma$ and is constructed of the given path. It allows us to define the integral for the irregular paths as an approximation by the standard line integrals over line segments that are patched together. In the next section we develop the integration theory over irregular paths for time- and path-dependent integrands in $\beta$-Hölder classes with $1 \geq \beta > 0$. To completely describe the set of irregular paths where the integral acts within a $\beta$-Hölder class, we introduce in addition to the $\alpha$-Hölder property of the path, $1 \geq \alpha \geq \beta$, a new independent characteristic of continuous functions based on Lévy's area between subsequent approximations of path $\gamma$, see Theorem 2.4 and Example 4.1. Similarly to the Gubinelli integration theory where for each $\alpha$-Hölder rough path the integral is defined on some subset of $\beta$-Hölder functions, here, in the opposite way, we show that to prove the existence of introduced integral for all $\beta$-Hölder integrands when $\beta \leq 1/2$, the set of irregular paths should necessarily be a proper subset of $\alpha$-Hölder class, see examples in Section 4 and Corollary 2.2. We also prove an important formula for applications that represents the integral over irregular paths as standard integrals, allowing the use of power of the existing fundamental theorem of calculus and table of indefinite integrals to evaluate it. With these results, in Section 3 we study the existence, uniqueness and continuity of solution of time- and path-dependent ODEIP in traditional $\beta$-Hölder classes, $0 < \beta \leq 1$. These results are a supplement to the Lyons-Gubinelli integration theory, in particular, for Hölder exponents less than 1/3. In Section 4 we present necessary proofs.

## 2. One-dimensional integration theory

1. Integral existence. We define an integral

$$\int_a^b f(t, g(t)) \, dg(t), 0 \leq a < b \leq 1,$$

for deterministic functions $f, g$. Consider a contour $\gamma$, which is given by a parametric presentation $\gamma(t) = (t, g(t)), t \in [0,1]$, in the plane $(t, x)$ where $g$ is continuous on $[0,1]$ ($g \in \mathcal{C}([0,1] \to R)$) with the norm

$$|g, \mathcal{C}([0,1])| = \sup_{t \in [0,1]} |g(t)|,$$

and a differential form $\omega = f(t, x) dx$, where $f \in \mathcal{C}([0,1] \times R \to R)$. With these notations, the last expression can be written as

$$\int_a^b f(t, g(t)) \, dg(t) = \int_{\gamma([a.b])} f(t, x) dx, 0 \leq a < b \leq 1.$$

To define an integral over such a contour $\gamma$, we construct a special approximation sequence of segments in the plane that converges to $\gamma$, and for any continuous form $\omega$ over such a piecewise smooth curve, we use the standard line integral. The main element in the construction of a special approximation sequence of line segments is

**Definition 2.1.** Let $g: [0,1] \to R$ be a continuous function and $\gamma(t) = (t, g(t)), t \in [0,1]$, be a Jordan curve in the plane. Let $\mathcal{D}_k = \{t_{k,n} = n2^{-k}; n = 0, \ldots, 2^k\}$ and $\mathcal{D} = \bigcup_{k>0} \mathcal{D}_k$ be dyadic subdivisions of interval $[0,1]$,

$$h_{k,n} = h_{k,n}(g) = 2^k \int_{t_{k,n}}^{t_{k,n+1}} g(\tau) d\tau \quad (0 \le n < 2^k).$$

For $0 \le a < b \le 1, k = 1, 2, \ldots$, define

$$[n_k(a), n_k(b)] = \{n; \ 0 \le n \le 2^k, [t_{k-1,[n/2]}, t_{k-1,[n/2]+1}] \subset [a,b]\},$$

($[t_{k-1,[n/2]}, t_{k-1,[n/2]+1}]$ is the parent of $[t_{k,n}, t_{k,n+1}]$)

$$\gamma_{k,n} = [(t_{k,n}, h_{k,n}), (t_{k,n+1}, h_{k,n})] \cup [(t_{k,n+1}, h_{k,n}), (t_{k,n+1}, h_{k,n+1})], n < n_k(b),$$
$$\gamma_{k,n} = [(t_{k,n}, h_{k,n}), (t_{k,n+1}, h_{k,n})], n = n_k(b),$$

and

$$\gamma_k = \gamma_k(a,b,g) = \bigcup_{n=n_k(a)}^{n_k(b)} \gamma_{k,n},$$

with orientation from $(a, g(a))$ to $(b, g(b))$.

In the following definition of the integral, we use the sequence $\{\gamma_k(g)\}$ as an approximator for the curve $\gamma$.

**Definition 2.2.** For $0 \le a < b \le 1$, consider the sequence

$$\int_{\gamma_k(a,b,g)} f(t,x) dx, \ k = 1,2,\ldots,$$

where $f \in C([0,1] \times R)$, $g \in C([0,1])$, $\gamma_k(a,b,g)$ are defined in Definition 2.1, and the integral over $\gamma_k(a,b,g)$ is a standard line Riemann integral. If a limit of this sequence exists, then we define the integral

$$\int_a^b f(t, g(t)) \, dg(t) = \lim_{k\to\infty} \int_{\gamma_k(a,b,g)} f(t,x) dx = \lim_{k\to\infty} \sum_{n=n_k(a)+1}^{n_k(b)} \int_{h_{k,n-1}}^{h_{k,n}} f(t_{k,n}, x) dx. \tag{2.1}$$

Integral (2.1) is properly defined since it is the limit of line Riemann integrals.
Immediately from the definition we can see that if the function $f$ depends only on $x$ then the integral exists for any continuous $f, g$ and

$$\int_a^b f(g(t)) \, dg(t) = \lim_{k\to\infty} \sum_{n=n_k(a)+1}^{n_k(b)} \int_{h_{k,n-1}}^{h_{k,n}} f(x) dx = \int_{g(a)}^{g(b)} f(x) dx,$$

see also (2.10) after Theorem 2.5 later. It is also obvious that this case can be extended to summable functions $f$ and continuous $g$.

**Theorem 2.1.a.** Let the function $f(t,x) = f(x)$ then integral (2.1) exists for any continuous function $g$ and any function $f$ summable on every compact subset of $R$.

In this setting, we mainly work with Hölder functions on a compact $K \subset [0,1]$ with the seminorm

$$|f, H_\beta(K)| = \sup_{x,y \in K, x \ne y} \frac{|f(x) - f(y)|}{|x-y|^\beta} < \infty, 0 < \beta \le 1.$$

The norm that we use in the space of the Hölder functions is given by $|f(x_0)| + |f, H_\beta(K)|, x_0 \in K$.
In the next theorem, we provide the condition for the existence of integral (2.1) in terms of a metric characteristic of the function $g$.

**Theorem 2.1.b.** Let the function $f$ be from Definition 2.2 and $f(\cdot, x) \in H_\beta, \beta > 0$, locally on $R$ uniformly by $x$,

$$|f, CH_\beta(K_g)| = \sup_{x \in [c,d]} |f(\cdot, x), H_\beta([0,1])| < \infty, g(t) \in [c,d], t \in [0,1], K_g = [0,1] \times [c,d].$$

Then integral (2.1) exists for any $g \in H_\alpha, \alpha > 0$, if

$$\sum_{k=1}^{\infty} 2^{k(1-\beta)} |B_{k-1}(g)| < \infty, \qquad (2.2)$$

where $B_k(g)$ is Levy's area between curves $\gamma_k(g)$ and $\gamma_{k+1}(g)$. Also

$$\left| \int_a^b f(t, g(t)) dg(t) - \int_{\gamma_{k_0+1}(a,b,g)} f(t, x) dx \right|$$

$$\leq |f, CH_\beta(K_g)| \sum_{k=k_0+1}^{\infty} 2^{k(1-\beta)} |B_{k-1}(g)| + 8|f, C(K_g)| \, |g, H_\alpha| \, |b - a|^\alpha, \qquad (2.3)$$

where $0 \leq a < b \leq 1, 2^{-k_0-1} \leq (b - a) \leq 2^{-k_0+1}$, and

$$|f, C(K_g)| = \sup_{(t,x) \in K_g} |f(t, x)|.$$

The construction of the sequence $\{B_k(g)\}$ is similar to the boundary stratification introduced by Salaev et al. (1990) to solve the problem on invariance of the singular integral of Cauchy type in Hölder classes.

For $g \in H_\alpha([0,1])$, $|B_k(g)| \leq C 2^{-k\alpha}$; therefore, for Hölder functions, the last theorem repeats the Young (1936) result for the new integral.

**Corollary 2.1**. With conditions from the previous theorem, if $\beta + \alpha > 1$, then integral (2.1) exists. For W-a.a. BMSPs integral (2.1) exists if $\beta > 1/2$.

Theorem 2.1.b states that condition (2.2) on the function $g$ is sufficient for the existence of integral (2.1) for any $\beta > 0$. In the following, we show that condition (2.2) is also necessary. We do so by studying integral (2.1) when the function $f$ depends only on the first argument

$$\int_a^b f(t) \, dg(t). \qquad (2.4)$$

**Theorem 2.2.** Let $g, \gamma$ be the same as in Definition 2.1, $0 < \beta < 1$ and $g \in H_\alpha([0,1]), 0 < \alpha < 1$; then, integral (2.4) exists for any function $f \in H_\beta([0,1])$ if and only if condition (2.2) holds.

The sufficiency part of Theorem 2.2 follows from Theorem 2.1.b. In the proof of necessity, for any $g \in C$, we construct a function $f \in H_\beta([0,1])$ such that for $k_0 > 1$,

$$\int_{\gamma_{k_0}(0,1,g)} f(t) dg(t) = \sum_{k=1}^{k_0-1} 2^{(k+1)(1-\beta)} |B_k(g)|,$$

which, with the assumption that integral (2.4) exists for every function from $H_\beta([0,1])$, implies condition (2.2). From theorems 2.1.b and 2.2, it immediately follows that condition (2.2) is also necessary and sufficient for the existence of integral (2.1). For $\beta > 0$, define the sets

$$GB_{\beta,\alpha}([0,1]) = \{g \in H_\alpha; \, \forall f \in H_\beta \, \exists \int_0^1 f(t) dg(t)\};$$

$$GB_{\beta,\alpha} = \{g \in H_\alpha; \, \forall f \, from \, Theorem \, 2.1.b \, \forall a \, \forall b, 0 \leq a < b \leq 1, \exists \int_a^b f(t, g(t)) dg(t)\};$$

$$GH_\beta = \{g \in C; \, \sum_{k=1}^{\infty} 2^{k(1-\beta)} |B_{k-1}(g)| < \infty\}.$$

With these notations,

$$GH_\beta \cap H_\alpha \subset GB_{\beta,\alpha} \subset GB_{\beta,\alpha}([0,1]) \subset GH_\beta \cap H_\alpha,$$

and all sets are equal, which means that condition (2.2) is a complete metric description of the existence of integral (2.1) for functions $f \in H_\beta$ and $g \in H_\alpha$.

**Theorem 2.2'.** Let $g \in H_\alpha$, the function $f$ be from Theorem 2.1.b, and $0 < \alpha, \beta < 1$. Then, integral (2.1) exists for any function $f$, $\sup_{x \in [c,d]} |f(\cdot, x), H_\beta([0,1])| < \infty$, if and only if condition (2.2) holds.

2. Integral calculus. The main properties of integral (2.1) are its continuity by the upper limit and the tools to facilitate its calculation.

**Theorem 2.3.** For $g \in GH_\beta \cap H_\alpha, 0 < \alpha, \beta < 1$, integral (2.1) exists on each interval $[a,b] \subset [0,1]$ for any function $f$, $\sup_{x \in [c,d]} |f(\cdot, x), H_\beta([0,1])| < \infty$, (see notations in the statement of Theorem 2.1.b) and the map

$$t \to \int_0^t f(s, g(s)) dg(s) \tag{2.5}$$

is continuous on the interval $[0,1]$. For $\beta > 2^{-1}$ map (2.5) is continuous for W-a.a. BMSPs.

The existence of integral (2.1) on each interval $[a,b] \subset [0,1]$ follows from Theorem 2.1. With condition (2.2), map (2.5) is continuous on the interval $[0,1]$ that follows from (2.3) in Theorem 2.1.b and the equality

$$\int_a^b f(s, g(s)) dg(s) = \int_a^c f(s, g(s)) dg(s) + \int_c^b f(s, g(s)) dg(s), \tag{2.6}$$

$0 \leq a < c < b \leq 1$, which will be proven in Section 4.

In Section 3 we use operator which relates to map (2.5)

$$(T_x y)(t) = y_0 + \int_0^t y(u) dx(u), t \in [0, T] \subset [0,1], y \in H_\beta, y_0 \in R, \tag{2.7}$$

with more complex form given by a function $F(u, y, x)$. The main properties of $T_x$ that we use studying ODEIP are shown in

**Theorem 2.4.** Let $\alpha, \beta > 0$. Then, $T_x: H_\beta \to H_\beta$ for any $x \in H_\alpha$ if and only if $\alpha \geq \beta$ and

$$\sum_{k=k_0+1}^{\infty} 2^{k(1-\beta)} |B_{k-1}(x)| \leq C 2^{-\beta k_0}, \tag{2.8}$$

where $C$ does not depend on $k_0 > 1$.

The sufficiency part follows from Theorem 2.1.b. If $T_x: H_\beta \to H_\beta$ then from Theorem 2.2 for constructed $f \in H_\beta$

$$C(b-a)^\beta \geq |(T_x f)(b) - (T_x f)(a)| \geq \sum_{k=k_0+1}^{\infty} 2^{k(1-\beta)} |B_{k-1}(x)|,$$

where $k_0$ is such that $2^{-k_0-1} \leq (b-a) \leq 2^{-k_0+1}$ which proofs (2.8). Similarly, with $y \equiv 1$ on $[0,1]$ and for any $0 \leq a < b \leq 1$

$$C |b-a|^\beta \geq |(T_x y)(b) - (T_x y)(a)| = |\int_a^b dx(u)| = |x(a) - x(b)|,$$

which implies that any $x \in H_\alpha$ is also from $H_\beta$, thus, $\alpha \geq \beta$.

The fundamental theorem of calculus is the basic computational tool in one-variable analysis. In two dimensions, such a tool is the Green formula (4.6). We prove the Green formula for special domains (Lemma 4.3) to reduce the calculation of integral (2.5) to Riemann (or Lebesgue) integrals in

**Theorem 2.5.** Let $\gamma$ be a contour, given by a parametric presentation $\gamma(t) = (t, g(t)), t \in [0,1]$ in the plane, where $g \in H_\alpha, g(0) = 0$, $\theta = \theta(s)$ is a segment, connecting points $\gamma(0)$ and $\gamma(s) = (s, g(s))$, $s \in (0,1]$. This segment consists of a closed set $\theta \cap \gamma$ and open segments $J_k = (\gamma(a_k), \gamma(b_k)), k > 0, \cup_k J_k = \theta \setminus \gamma$, such that $(+\gamma([a_k, b_k])) \cup (-J_k)$ or $(-\gamma([a_k, b_k])) \cup (+J_k)$ is the Jordan boundary of a bounded domain $G_k$. Let the function $f$ be the same as in Theorem 2.1.b, $f \in ACL$ (Absolutely Continuous on Lines), and $\partial f / \partial t \in L(K_g)$. Then, for $g \in GH_\beta \cap H_\alpha, \alpha, \beta > 0$,

$$\int_0^s f(t, g(t)) dg(t) = \sum_k \sigma_k \int_{G_k} \frac{\partial}{\partial t} f(t, x) dt dx + \frac{g(s)}{s} \int_0^s f\left(t, \frac{g(s)}{s} t\right) dt, \tag{2.9}$$

where $\sigma_k = 1$ if $\partial G_k = (+\gamma([a_k, b_k])) \cup (-J_k)$, and $\sigma_k = -1$ if $\partial G_k = (-\gamma([a_k, b_k])) \cup (+J_k)$. For $\beta > 2^{-1}$, formula (2.9) is valid for W-a.a. BMSPs.

For functions $f(t, x) = f(x)$ formula (2.9) has a simple expression

$$\int_0^s f(g(t))dg(t) = \int_0^{g(s)} f(t)dt, \qquad (2.10)$$

or for any continuous function $g$ and differentiable function $F$ such that $F'$ exists a.e. on $R$ and is locally summable

$$F(g(s)) - F(g(0)) = \int_0^s \frac{d}{dx}F(g(t))dg(t),$$

which resembles stochastic 1-dimensional Itô's formula for the Stratonovich integral. Since integral (2.10) exists for any continuous function $g$, we can have a.e. on $\Omega$ an integral

$$\int_a^b f(X(t,\omega))dX(t,\omega)$$

for any stochastic process $X(t,\omega)$ with continuous sample paths on a probability space $(\Omega, \mathcal{F}, P)$, particularly for Brownian motion or fractional Brownian motion with Hurst parameter $H \in (0,1)$. In the theorem below we compare this integral with the Itô and Stratonovich stochastic integrals.

**Theorem 2.6.** Let $(\Omega, \mathcal{F}, P)$ be a probability space, $X(t,\omega)$ be a continuous process on it, and $f$ be continuously differentiable on $R$. Then the integral

$$\int_0^s f(X(t,\omega))dX(t,\omega)$$

exists for a.a. $\omega \in \Omega$, is continuous on $[0,1]$, and

$$\int_0^s f(X(t,\omega))dX(t,\omega) = ((\text{Itô}) \int_0^s f(X(\tau))dX(\tau))(\omega) + \frac{1}{2}\int_0^s \frac{df}{dx}(X(\tau))d\tau,$$

$$\int_0^s f(X(t,\omega))dX(t,\omega) = ((\text{Strat}) \int_0^s f(X(\tau))dX(\tau))(\omega),$$

for a.a. $\omega \in \Omega$ and

$$(\text{Itô/Strat}) \int_0^s f(X(\tau))dX(\tau)$$

are continuous versions of corresponding stochastic integrals if they exist.

Proof. Let

$$g(x) = \int_0^x f(\tau)d\tau, Y(t) = g(X(t)).$$

By Itô's formula,

$$g(X(s)) - g(X(0)) = (\text{Itô})\int_0^s f(X(\tau))dX(\tau) + \frac{1}{2}\int_0^s \frac{df}{dx}(X(\tau))d\tau,$$

and by formula (2.10),

$$g(X(s,\omega)) - g(X(0,\omega)) = \int_0^s f(X(t,\omega))dX(t,\omega),$$

which gives us the first equality from the theorem. The second one follows from the last and the relationship between Itô and Stratonovich integrals:

$$(\text{Strat})\int_0^s f(X(\tau))dX(\tau) = (\text{Itô})\int_0^s f(X(\tau))dX(\tau) + \frac{1}{2}\int_0^s \frac{df}{dx}(X(\tau))d\tau.$$

For $f(t,x) = f(t)$, formula (2.9) is the integration by parts

$$f(s) * g(s) = \int_0^s f(t)dg(t) + \int_0^s g(t)df(t), s \in [0,1]. \qquad (2.11)$$

These formulas allow the use of all powers of tables of indefinite integrals to calculate integrals over irregular paths. For example,

$$\int_0^s \sin(a\,x(t))e^{b\,x(t)}dx(t) = \frac{e^{b\,x(s)}}{a^2+b^2}(b\sin(a\,x(s)) - a\cos(a\,x(s))), a^2+b^2 > 0,$$

$$\int_0^s |x(t)|^\beta dx(t) = \frac{1}{\beta+1}\text{sign}(x(s))|x(s)|^{\beta+1},$$

and the famous one

$$\int_0^s x(t)dx(t) = \frac{1}{2}(x(s))^2.$$

The stochastic Itô and Stratonovich formulas mentioned in Theorem 2.6 were initially proved for the Itô (1951) and Stratonovich (1966) stochastic integrals. Both of these integrals have their own advantages and are widely used in many applications. For the Itô integral, the formula, named Itô's formula, has an additional term depending on the second derivative that the Stratonovich and the new integrals (formula (2.10)) do not have.

The calculation formula (2.9) gives another presentation of the new integral as standard Riemann or Lebesgue integrals, allowing the use of power of the existing fundamental theorem of calculus and table of indefinite integrals to evaluate the new integral. For example, we can use calculation formulas (2.9) – (2.11) to obtain the (numerical) solutions of ODEIP like (3.2) as ordinary integrals and analyze them with standard calculus techniques.

3. Integration in $\alpha$-Hölder classes with $\alpha \leq 2^{-1}$. Properties of operator $T_x$ from (2.7) play essential role in the study of ODEIP. From Theorem 2.4 we see that conditions $\alpha \geq \beta$ and (2.8) are complete for $T_x$ to act from $H_\beta$ into $H_\beta$ for any $x \in H_\alpha$. Also, if $\alpha + \beta > 1$ then (2.8) holds and again $T_x: H_\beta \to H_\beta$, thus, the only condition $\alpha + \beta > 1$ allows to conclude that $T_x$ acts in $H_\beta$. As shown in Example 4.2 this is not true if $\alpha + \beta < 1$: for any $\alpha, \beta > 0$ and $x \in H_\alpha$ there is $y \in H_\beta$ such that $T_x y$ does not exist, if $\alpha + \beta < 1$. The last case is very important due to applications to BMSPs. Wiener (1930) proved that W-a.a BMSPs are $\alpha$-Hölder with $\alpha < 2^{-1}$ and the next theorem shows that if $\beta \leq 2^{-1}$ then integral (2.1) over W-a.a. BMSPs cannot exist for all $f \in H_\beta$. The proof is similar to that in Cameron and Martin (1947) and presented in Section 4.

**Theorem 2.7.** For W-a.a. $g \in C$,

$$\lim_{k\to\infty} 2^{-k/2} \sum_{m=0}^{2^k-1} |h_{k+1,2m+1}(g) - h_{k+1,2m}(g)| = \sqrt{\frac{2}{3\pi}},$$

where $h_{k,m}$ is from Definition 2.1.
Since

$$GH_\beta = \{g \in C; \sum_{k=1}^\infty 2^{-k\beta} \sum_{m=0}^{2^k-1} |h_{k+1,2m+1}(g) - h_{k+1,2m}(g)| < \infty\},$$

then $W(GH_\beta) = 1$ for $\beta > 2^{-1}$ and $W(GH_\beta) = 0$ for $\beta \leq 2^{-1}$. Thus,

integral (2.1) exists for W-a,a, BMSPs and all functions from $H_\beta$ if and only if $\beta > 2^{-1}$, and for W-a.a. BMSPs there exists a function from $H_\beta$ such that integral (2.1) does not exist if and only if $\beta \leq 2^{-1}$.

This theorem completes Young's result for Hölder functions in the general case $f(t,x)$ and, in particular, when $f$ does not depend on $x$, which relates to the Riemann-Stieltjes integration.

**Corollary 2.2.** Integrals (2.1) and, consequently, (2.4) exist for W-a.a. BMSPs and for any function from $H_\beta$ if and only if $\beta > 2^{-1}$ and for W-a,a BMSPs there exists a function from $H_\beta$ such that integral (2.1) does not exist if and only if $\beta \leq 2^{-1}$.

Thus, if one wants to develop an integration theory over irregular paths from the entire Hölder class $H_\alpha$ for $\beta$-Hölder functions with $\alpha, \beta \leq 2^{-1}$, a proper subclass of $H_\beta$ has to be constructed. Such approach was developed in the pioneering work of Lyons (1998), (2002) for differential equations driven by $\alpha$ −Hölder rough paths, $\alpha > 1/3$, where a technique was developed based on a given values of higher indefinite iterated integrals and a topology was constructed to define such subspace as a closure of smooth functions in this topology. More explicitly, Gubinelli (2004) introduced a subclass of $\beta$−Hölder functions, $\beta > 1/3$, whose increments are controlled by a Lyons rough path where an integral can be defined and differential equations driven by such paths can be studied.

In this paper we developed an integration theory that allowed us to describe all continuous irregular paths for which the new integral is defined for all functions from $H_\beta, \beta > 0$. As we see from Corollary 2.2, this integral over BMSPs does not exist for all functions from $H_\beta$ when $0 < \beta \leq 1/2$. But there are many other irregular paths in the class $H_\alpha, 0 < \alpha \leq 2^{-1}$, like in Example 4.1, where the new integral is defined for all functions from $H_\beta, \alpha + \beta \leq 1$, and this integral will be applied in the next section to study ODEIP over such irregular inputs.

For $\beta > 2^{-1}$, integral (2.1) can be used over W-a.a. BMSPs and for all functions $f \in H_\beta$. If there is still a need for some $g \in \mathcal{C}$ to verify the applicability of the integral, condition (2.2) has to be checked.

For $\beta \leq 2^{-1}$, integral (2.1) cannot be used over W-a.a. BMSPs for all functions $f \in H_\beta$. However, again, if there is still a need for some $g \in \mathcal{C}$ to verify the applicability of the integral, condition (2.2) has to be checked.

Thus, condition (2.2) which is also necessary for the existence of integral (2.1), is the most useful tool to verify the applicability of integral (2.1) for Hölder classes.

In the next theorem we show that in the constructed Gubinelli type classes the new integral exists for all paths from the given class $H_\alpha, \alpha > 0$.

We will use generalized Hölder classes $H_{\psi(\delta)}$ on a compact $K \subset [0,1]$ with the seminorm

$$|f, H_{\psi(\delta)}(K)| = \sup_{x,y \in K, x \neq y} \frac{|f(x) - f(y)|}{\psi(|x-y|)} < \infty,$$

where $\psi: (0,1] \to R_+$ such that $\psi$ is nondecreasing and $\psi(\delta) \to 0$ as $\delta \to 0$. For $\psi(\delta) = \delta^\alpha$ this is a traditional Hölder class $H_\alpha$.

**Definition 2.3.** (Gubinelli classes). A function $y \in \mathcal{C}([0,1])$ is said to be controlled by $x \in \mathcal{C}([0,1])$ if there exists a function $y' \in H_{\psi(\delta)}([0,1])$ such that

$$|y(t_{k+1,2n}) - y(t_{k+1,2n+1}) - y'(t_{k,n})(h_{k+1,2n}(x) - h_{k+1,2n+1}(x))|$$
$$\leq C |y', H_{\psi(\delta)}([0,1])|\psi(2^{-k})|h_{k+1,2n}(x) - h_{k+1,2n+1}(x)|,$$
$$|y(t_{k+1,2n+1}) - y(t_{k+1,2(n+1)}) - y'(t_{k+1,2n+1})(h_{k+1,2n}(x) - h_{k+1,2n+1}(x))|$$
$$\leq C |y', H_{\psi(\delta)}([0,1])|\psi(2^{-k})|h_{k+1,2n}(x) - h_{k+1,2n+1}(x)|,$$

$k > 1, n = 0, \ldots, 2^k - 1, t_{k,n} = n2^k, h_{k,n}(x)$ are defined in Definition 2.1 with $x = g$ and the constant $C$ does not depend on $y$ and $y'$.

The set of all such functions $y$ we denote as $\mathbb{G}_{x,y'}([0,1])$.

**Theorem 2.8.** Let $f \in \mathbb{G}_{g,f'}, f' \in H_{\psi(\delta)}, g \in H_{\delta^\alpha}, \alpha > 0$. Then integral (2.4) exists if

$$\sum_{k=1}^\infty \psi(2^{-k}) \sum_{n, [t_{k,n-1}, t_{k,n}] \subset [a,b]} \left(h_{k+1,2(n-1)}(g) - h_{k+1,2(n-1)+1}(g)\right)^2 < \infty.$$

**Corollary 2.3.** For W-a.a. $g \in \mathcal{C}$

$$\sum_{n, [t_{k,n-1}, t_{k,n}] \subset [a,b]} \left(h_{k+1,2(n-1)}(g) - h_{k+1,2(n-1)+1}(g)\right)^2 \leq K < \infty, k > 1,$$

see Cameron and Martin (1947) and Theorem 2.7 above, and integral (2.4) exists if

$$\sum_{k=1}^\infty \psi(2^{-k}) < \infty. \tag{2.12}$$

**Corollary 2.4.** If $g, f' \in H_{\delta^\alpha}$, $\alpha > 1/3$, then integral (2.4) exists for any $f \in \mathbb{G}_{g,f'}$.

From Definition 2.3 and the Marchaud (1927) theorem, for $f', g \in H_\alpha$, $\mathbb{G}_{g,f'} \subset H_{\delta^{2\alpha}}$. This means that for $f \in \mathbb{G}_{g,f'} \subset H_{2\alpha}, g \in H_\alpha$, and $\alpha > 1/3$, since $2\alpha + \alpha > 1$, in the proof of Corollary 2.4 we can use the classical Young integration theory for Hölder classes.

In general, it would be interesting to show a class $\mathbb{G}_{g,f'}$ where integral (2.4) exists, and it is not included in any $H_{\delta^\gamma}$ with $\alpha + \gamma > 1$ ($g \in H_\alpha$). Thus, for $f' \in H_{\psi(\delta)}$, $\psi(\delta) = (\ln 1/\delta)^{-1}(\ln \ln 1/\delta)^{-2}$, the condition (2.12) holds, and integral (2.4) exists for W-a.a. BMSPs $g$ and any function from $\mathbb{G}_{g,f'}$. Moreover, with Lévy's modulus of continuity theorem, Lévy (1937), and by the Marchaud theorem, $\mathbb{G}_{g,f'} \subset H_{\theta(\delta)}$, where $\theta(\delta) = \delta^{1/2}(\ln \ln 1/\delta)^{-2}$, and from Theorem 2.7

$$\sum_{k=1}^{k_0-1} \theta(2^{-k}) \sum_{m=0}^{2^k-1} |h_{k+1,2m+1}(g) - h_{k+1,2m}(g)| \geq C \sum_{k=1}^{k_0-1} (\ln(k \ln 2))^{-2} \to \infty, \text{ as } k_0 \to \infty.$$

Having this, based on the proof of Theorem 2.2, for W-a.a. BMSPs integral (2.4) does not exists for all functions from $H_{\theta(\delta)}$. Also, the minimum class $H_{\delta^\beta}$ that includes $H_{\theta(\delta)}$ is $H_{\delta^{1/2}}$ and the maximum class $H_{\delta^\alpha}$ that BMSPs includes is also $H_{\delta^{1/2}}$ that gives us $\alpha + \beta = \frac{1}{2} + \frac{1}{2} = 1$. For such case we cannot guarantee applicability of the Young integration theory and should work with the class $\mathbb{G}_{g,f'}$.

## 3. Ordinary differential equations over irregular paths

We use the results from the previous section to apply the introduced integral to prove the existence, uniqueness and continuity of solution of ODEIP in traditional Hölder classes. Consider $F: R_+ \times R^m \times R^d \to \mathcal{L}(R^d, R^m), x: [0, T] \to R^d, y: [0, T] \to R^m$, and differential equations

$$dy(t) = F(t, y(t), x(t))dx(t), y(0) = y_0, t \in [0, T]. \quad (3.1)$$

The solutions of this equations we will study in Hölder spaces with the norm

$$|y_0| + |y, H_\beta([0,T])| = |y_0| + \sup_{t,u \in [0,T], t \neq u} \frac{|y(t) - y(u)|}{|t-u|^\beta},$$

with $x$ from the space of continuous functions

$$G_{\beta,\gamma} = \{x \in C([0,T] \to R^d); |x, G_{\beta,\gamma}| < \infty\}, \beta, \gamma > 0,$$

$$|x, G_{\beta,\gamma}([0,T])| = \sup_{\substack{0 \leq a < b \leq T, \\ 1 \leq j \leq d}} \frac{1}{(b-a)^{\beta+\gamma}} \mu(x_j, a, b, \beta),$$

$$\mu(g, a, b, \beta) = \sum_{k=k_0+1}^\infty 2^{-(k+1)\beta} \sum_{n, [t_{k,n-1}, t_{k,n}] \subset [a,b]} |h_{k+1,2(n-1)}(g) - h_{k+1,2(n-1)+1}(g)| < \infty,$$

where $\mu(x_j, a, b, \beta)$ is a functional characteristic of the curve from Estimate 4.1 for components of vector $x$, $2^{-k_0-1} \leq (b-a) \leq 2^{-k_0+1}$.

In Example 4.1 later we show that the set $\{x \in C([0,T] \to R^d); |x, G_{\beta,\gamma}([0,T])| < \infty\}$ includes highly oscillatory irregular path.

Below, we will use also the obvious bound

$$|x, G_{\beta_1,\gamma}([0,T])| \leq 2|x, G_{\beta_2,\gamma}([0,T])|, 0 < \beta_2 < \beta_1.$$

For $x \in G_{\beta,\gamma} \cap H_\alpha, \beta \leq \alpha$, any function $y \in H_\beta$ that satisfies the integral equation

$$y(t) = y_0 + \int_0^t F(u, y(u), x(u))dx(u), t \in [0, T], \quad (3.2)$$

if it exists, is a solution of differential equation (3.1). Here the notation

$$\int_a^b F(u, y(u), x(u))dx(u)$$

denotes an m-dimensional vector with $i$-entry, $1 \leq i \leq m$,

$$\int_a^b F_i(u, y(u), x(u))dx(u) = \sum_{j=1}^d \int_a^b F_{ij}(u, y(u), x(u))dx_j(u),$$

and the last integral is understood in the sense of Definition 2.2, and its existence was proved in Theorem 2.1 for $f(t, x_j) = F_{ij}(t, y(t), x_1(t), ,, x_j ... x_d(t))$.

We follow Gubinelli (2004) to prove existence, uniqueness and continuity of solutions of ODEIP in Hölder classes based on the new integral.

**Lemma 3.1**. Let the function $F$ be $F: R_+ \times R^m \times R^d \to \mathcal{L}(R^d, R^m)$, its components $F_{ij} \in \mathcal{C}(R_+ \times R^m \times R^d)$ and $x \in H_\alpha([0,T], R^d)$, $|x, \mathcal{C}([0,T])| \le M < \infty$, $y, \tilde{y} \in H_\beta([0,T], R^m)$, $|y, \mathcal{C}([0,T])| \le M$, $|\tilde{y}, \mathcal{C}([0,T])| \le M < \infty$, $\beta \le \alpha$.

a) If for all $t, y_1, \ldots, y_m, x_1, \ldots, \hat{x}_j, \ldots, x_d$ ($x_j$ is excluded) $F_{ij}$ is $\theta$-Hölder, $0 < \theta \le 1$, locally uniformly by $x_j$ then the function $f(t, x_j) = F_{ij}(t, y(t), x_1(t), \ldots, x_j, \ldots, x_d(t))$ is $\theta\beta$-Hölder by $t$ for all $x_j$ and
$$\sup_{x_j}|f(\cdot, x_j), H_{\theta\beta}| \le |F_{ij}, H_\theta(t)| T^{\theta(1-\beta)} + |F_{ij}, H_\theta(y)||y, H_\beta| + |F_{ij}, H_\theta(x)||x, H_\alpha| T^{\theta(\alpha-\beta)} < \infty,$$
where $y = (y_1, \ldots, y_m), x = (x_1, \ldots, x_d)$,
$$|F_{ij}, H_\theta(y)| = \sum_{k=1}^m |F_{ij}, H_\theta(y_k)|, \quad |F_{ij}, H_\theta(x)| = \sum_{k=1,,k \neq j}^d |F_{ij}, H_\theta(x_k)|,$$
$$|F_{ij}, H_\theta(t)| = \sup_{|x| \le M, |y| \le M} |F_{ij}(\cdot, y, x), H_\theta|, \quad |F_{ij}, H_\theta(y_k)| = \sup_{t \in [0,1], |x| \le M, |y| \le M} |F_{ij}(t, \ldots, y_{k-1}, \cdot, \ldots, x), H_\theta|,$$
$$|F_{ij}, H_\theta(x_k)| = \sup_{t \in [0,1], |x| \le M, |y| \le M} |F_{ij}(t, y, \ldots, x_{k-1}, \cdot, \ldots), H_\theta|, k \neq j.$$

b) With conditions from section a), let $0 < \sigma < 1$, then
$$|F_{ij}(\cdot, y(\cdot), x(\cdot, \hat{x}_j)) - F_{ij}(\cdot, \tilde{y}(\cdot), x(\cdot, \hat{x}_j)), H_{\sigma\beta\theta}([0,T])|$$
$$\le (|F_{ij}, H_\theta(y)|(|y, H_\beta| + |\tilde{y}, H_\beta|) + 2(|F_{ij}, H_\theta(t)| T^{\theta(1-\beta)} + |F_{ij}, H_\theta(x)||x, H_\alpha| T^{\theta(\alpha-\beta)}))^\sigma$$
$$* (2|F_{ij}, H_\theta(y)||y - \tilde{y}, \mathcal{C}|^\theta)^{1-\sigma}.$$

c) For $\theta$-Hölder continuously differentiable by $y$ function $F_{ij}(t, y, x)$,
$$|F_{ij}(\cdot, y(\cdot), x(\cdot, \hat{x}_j)) - F_{ij}(\cdot, \tilde{y}(\cdot), x(\cdot, \hat{x}_j)), H_{\beta\theta}([0,T])|$$
$$\le |D_y F_{ij}, \mathcal{C}| T^{\beta(1-\theta)}|y - \tilde{y}, H_\beta| + (|D_y F_{ij}, H_\theta(y)|(|y, H_\beta| + |\tilde{y}, H_\beta|)$$
$$+ |D_y F_{ij}, H_\theta(t)| T^{\theta(1-\beta)} + |D_y F_{ij}, H_\theta(x)||x, H_\alpha| T^{\theta(\alpha-\beta)})|y - \tilde{y}, \mathcal{C}|,$$
where $D_y F_{ij}: R_+ \times R^m \times R^d \to \mathcal{L}(R^m, R)$,
$$D_y F_{ij}(t, y, x) = (\frac{\partial}{\partial y_k} F_{ij}(u, y, x))_{k=1}^m.$$

Proof. a). By the lemma assumptions
$$|f(t_1, x_j) - f(t_2, x_j)|$$
$$\le |F_{ij}, H_\theta(t)||t_1 - t_2|^\theta + \sum_{k=1}^m |F_{ij}, H_\theta(y_k)||y_k, H_\beta([0,T])||t_1 - t_2|^{\theta\beta}$$
$$+ \sum_{k=1, k \neq j}^d |F_{ij}, H_\theta(x_k)||x, H_\alpha([0,T])||t_1 - t_2|^{\theta\alpha}.$$

b) Applying the interpolation arguments from Proposition 5, Gubinelli (2004), we first estimate in $\mathcal{C}$ norm
$$A(y, \tilde{y}, u, s) = |F_{ij}(u, y(u), x_1(u), \ldots, x_j, \ldots, x_d(u)) - F_{ij}(u, \tilde{y}(u), x_1(u), \ldots, x_j, \ldots, x_d(u))$$
$$- (F_{ij}(s, y(s), x_1(s), \ldots, x_j, \ldots, x_d(s)) - F_{ij}(s, \tilde{y}(s), x_1(s), \ldots, x_j, \ldots, x_d(s)))|.$$
Since the value of $x_j$ is the same we will omit $x_j$ in further estimations. Thus,
$$A(y, \tilde{y}, u, s) \le |F_{ij}(u, y(u), x(u)) - F_{ij}(u, \tilde{y}(u), x(u))| + |(F_{ij}(s, y(s), x(s)) - F_{ij}(s, \tilde{y}(s), x(s)))|$$
$$\le 2|F_{ij}, H_\theta(y)||y - \tilde{y}, \mathcal{C}|^\theta.$$
The next is the estimate for $\theta\beta$-Hölder norm
$A(y, \tilde{y}, u, s)$
$$\le |F_{ij}(u, y(u), x(u)) - F_{ij}(u, y(s), x(u))| + |F_{ij}(u, \tilde{y}(u), x(u)) - F_{ij}(u, \tilde{y}(s), x(u))|$$
$$+ |F_{ij}(u, y(s), x(u)) - F_{ij}(s, y(s), x(s))| + |F_{ij}(u, \tilde{y}(s), x(u)) - F_{ij}(s, \tilde{y}(s), x(s))|$$
$$\le |F_{ij}, H_\theta(y)|(|y, H_\beta| + |\tilde{y}, H_\beta|) |u - s|^{\theta\beta}$$
$$+ 2(|F_{ij}, H_\theta(t)| T^{\theta(1-\beta)} + |F_{ij}, H_\theta(x)||x, H_\alpha(x)| T^{\theta(\alpha-\beta)})|u - s|^{\theta\beta}.$$
And the interpolation between $\mathcal{C}$ and $\theta\beta$-Hölder norms

$$\frac{|A(y,\tilde{y},u,s)|}{|u-s|^{\sigma\theta\beta}} = \left(\frac{|A(y,\tilde{y},u,s)|}{|u-s|^{\theta\beta}}\right)^{\sigma}|A(y,\tilde{y},u,s)|^{1-\sigma}$$
$$\leq (|F_{ij}, H_\theta(y)|(|y, H_\beta| + |\tilde{y}, H_\beta|) + 2(|F_{ij}, H_\theta(t)|T^{\theta(1-\beta)} + |F_{ij}, H_\theta(x)||x, H_\alpha(x)|T^{\theta(\alpha-\beta)})^\sigma$$
$$* (2|F_{ij}, H_\theta(y)||y - \tilde{y}, C|^\theta)^{1-\sigma}.$$

c) For continuously differentiable function $\varphi$ we are going to use the line integral formula in $R^m$
$$\varphi(x) - \varphi(y) = \int_0^1 D\varphi(tx + (1-t)y)dt * (x - y)$$
$$= \int_0^1 \sum_{k=1}^m \frac{\partial}{\partial y_k}\varphi(tx + (1-t)y)dt(x_k - y_k), x, y \in R^m.$$

Now we estimate $\beta\theta$-Hölder norm for
$$A(y, \tilde{y}, u, s) = |F_{ij}(u, y(u), x(u)) - F_{ij}(u, \tilde{y}(u), x(u)) - (F_{ij}(s, y(s), x(s)) - F_{ij}(s, \tilde{y}(s), x(s)))|$$
$$\leq |F_{ij}(u, y(u), x(u)) - F_{ij}(u, \tilde{y}(u), x(u)) - (F_{ij}(u, y(s), x(u)) - F_{ij}(u, \tilde{y}(s), x(u)))|$$
$$+ |F_{ij}(u, y(s), x(u)) - F_{ij}(u, \tilde{y}(s), x(u)) - (F_{ij}(s, y(s), x(s)) - F_{ij}(s, \tilde{y}(s), x(s)))|$$
$$= |\int_0^1 D_y F_{ij}(u, ty(u) + (1-t)\tilde{y}(u), x(u))dt * (y(u) - \tilde{y}(u) - (y(s) - \tilde{y}(s)))$$
$$+ \int_0^1 \left(D_y F_{ij}(u, ty(u) + (1-t)\tilde{y}(u), x(u)) - D_y F_{ij}(u, ty(s) + (1-t)\tilde{y}(s), x(u))\right)dt * (y(s) - \tilde{y}(s))|$$
$$+ |\int_0^1 \left(D_y F_{ij}(u, ty(s) + (1-t)\tilde{y}(s), x(u)) - D_y F_{ij}(s, ty(s) + (1-t)\tilde{y}(s), x(s))\right)dt * (y(s) - \tilde{y}(s))|$$
$$\leq |D_y F_{ij}, C||y - \tilde{y}, H_\beta||u - s|^\beta + |D_y F_{ij}, H_\theta(y)|(|y, H_\beta| + |\tilde{y}, H_\beta|)|u - s|^{\beta\theta}|y - \tilde{y}, C|$$
$$+ (|D_y F_{ij}, H_\theta(t)||u - s|^\theta + |D_y F_{ij}, H_\theta(x)||x, H_\alpha||u - s|^{\alpha\theta})|y - \tilde{y}, C|. \blacksquare$$

The properties of the indefinite integral proved in Theorem 2.3 and Lemma 3.1 allow one to repeat the standard solution theory for ordinary differential equations in classes $H_\beta$ for any $\beta > 0$.

**Theorem 3.1. (Existence).** Let $0 < \beta \leq \alpha \leq 1$, the function $F: R_+ \times R^m \times R^d \to \mathcal{L}(R^d, R^m)$ is the same as in Lemma 3.1 b) with $\theta$- Hölder components, $0 < \theta \leq 1$. If $0 < \beta - \gamma < \beta_0 < \theta\beta, 0 < \gamma$, then for any $x \in G_{\beta_0,\gamma} \cap H_\alpha([0,T], R^d)$ there exists a solution from $H_\beta([0,T], R^m)$ of equation (3.1).

Proof. Let $y \in H_\beta([0,T], R^m)$. By Lemma 3.1 a), the function $f_{ij}(t, x_j) = F_{ij}(t, y(t), x_1(t),,, x_j ... x_d(t))$ is $\theta\beta$-Hölder and by the theorem assumption, $\beta_0 < \theta\beta$, the condition to apply Theorem 2.1.b is fulfilled and we can define the map

$$(T_{x,F}y)(t) = y_0 + \int_0^t F(u, y(u), x(u))dx(u), t \in [0, T], \tag{3.3}$$

and also, the closed unit ball
$$B_T(\beta) = \{y \in H_\beta([0,T]); y(0) = y_0, |y, H_\beta([0,T])| \leq 1\}, \tag{3.4}$$
which is not empty.

For $0 \leq a < b \leq T$
$$|\int_a^b F(u, y(u), x(u))dx(u)| \leq \sum_{i=1}^m \sum_{j=1}^d |\int_a^b F_{ij}(u, y(u), x(u))dx_j(u)|$$
$$\leq \sum_{i=1}^m \sum_{j=1}^d (|F_{ij}, H_\theta(t)|T^{\theta(1-\beta)} + |F_{ij}, H_\theta(y)||y, H_\beta| + |F_{ij}, H_\theta(x)||x, H_\alpha|T^{\theta(\alpha-\beta)})\mu(x_j, a, b, \beta\theta)$$
$$+ 8|F_{ij}, C||x, H_\alpha([0,T])||b - a|^\alpha)$$
$$\leq \sum_{i=1}^m \sum_{j=1}^d (|F_{ij}, H_\theta(t)|T^{\theta(1-\beta)} + |F_{ij}, H_\theta(y)||y, H_\beta| + |F_{ij}, H_\theta(x)||x, H_\alpha|T^{\theta(\alpha-\beta)})$$
$$* |x, G_{\beta\theta,\gamma}(x, [0,T])||b - a|^{\beta\theta+\gamma} + 8|F_{ij}, C||x, H_\alpha([0,T])||b - a|^\alpha)$$
$$\leq C(1 + |y, H_\beta| + |x, H_\alpha|)|x, G_{\beta\theta,\gamma}(x, [0,T])||b - a|^{\beta\theta+\gamma} + 8|F_{ij}, C||x, H_\alpha([0,T])||b - a|^\alpha, \tag{3.5}$$
where $C$ depends only on $|x, C|, F, |y, C|$. It is clear from the previous estimate that for $y \in B_T(\beta)$

$$|T_{x,F}y, H_\beta([0,T])| \le C_1 T^{\gamma-\beta(1-\theta)} + C_2 |x, H_\alpha([0,T])| T^{\alpha-\beta},$$

where constants $C_1, C_2$ do not depend on $y \in B_T(\beta)$ and bounded by $T \in (0,1]$. Based on the last estimate, if $\alpha > \beta$, a small enough $T = T_1$ can be chosen such that $T_{x,F}: B_T(\beta) \to B_T(\beta)$.

If $\alpha = \beta$ we should use another well-known approach working with the Hölder class $H_{\alpha_1}$ when $\alpha_1 < \alpha$ and $\alpha_1 \theta > \beta_0$. In this case with $x \in G_{\beta_0,\gamma} \cap H_\alpha([0,T], R^d)$ we can consider $x \in G_{\beta_0,\gamma} \cap H_{\alpha_1}$ and for $y \in B_T(\alpha_1)$ based on estimate (3.5)

$$|T_{x,F}y, H_{\alpha_1}([0,T])| \le C_1 T^{\alpha_1\theta+\gamma-\alpha_1} + C_2|x, H_{\alpha_1}([0,T])| \le C_1 T^{\alpha_1\theta+\gamma-\alpha_1} + C_2|x, H_\alpha([0,T])| T^{\alpha-\alpha_1},$$

and again, we can choose a small enough $T = T_1$ such that $T_{x,F}: B_T(\alpha_1) \to B_T(\alpha_1)$ which proves the invariance of $T_{x,F}$ on $B_T(\alpha_1)$.

To prove existence in Hölder classes $H_\beta$ based on Leray–Schauder–Tychonoff theorem, we need to show that the map $T_{x,F}$ is continuous in the $\beta$-Hölder topology. Let $y, \tilde{y} \in B_T(\beta), \beta \le \alpha$, then for $\sigma = \frac{\beta_0}{\beta\theta} < 1$ applying Theorem 2.2 and Lemma 3.1 b) ($y, \tilde{y} \in H_\beta \subset H_{\sigma\beta}$)

$$|(T_{x,F}y)(a) - (T_{x,F}y)(b) - ((T_{x,F}\tilde{y})(a) - (T_{x,F}\tilde{y})(b))|$$

$$= |\int_a^b (F(u, y(u), x(u)) - F(u, \tilde{y}(u), x(u))) dx(u)|$$

$$\le \sum_{i=1}^m \sum_{j=1}^d ((|F_{ij}, H_\theta(y)|(|y, H_\beta| + |\tilde{y}, H_\beta|) + 2(|F_{ij}, H_\theta(t)|T^{\theta(1-\beta)} + |F_{ij}, H_\theta(x)||x, H_\alpha|T^{\theta(\alpha-\beta)})^\sigma)$$

$$* (2|F_{ij}, H_\theta(y)||y - \tilde{y}, C|^\theta)^{1-\sigma} |x, G_{\beta_0,\gamma}(x, [0,T])| |b-a|^{\gamma+\beta_0}$$

$$+ 8|F_{ij}, H_\theta(y)||y - \tilde{y}, C|^\theta |x, H_\alpha([0,T])| |b-a|^\alpha),$$

and

$$|T_{x,F}y - T_{x,F}\tilde{y}, H_\beta| \le C(|y - \tilde{y}, C|^{\theta(1-\sigma)} + |y - \tilde{y}, C|^\theta) \to 0, |y - \tilde{y}, H_\beta| \to 0,$$

where constant $C$ does not depend on $y, \tilde{y} \in B_T(\beta)$ and bounded by $T \in (0,1]$, which completes the proof of existence for $\beta < \alpha$.

For $\beta = \alpha$, since the proof of continuity of the map $T_{x,F}$ in the $\alpha_1$-Hölder topology for $\alpha_1 < \alpha$ is the same as the previous arguments, we have the existence of the solution of the equation (3.1) or (3.2) in classes $H_{\alpha_1}$. But if $y \in H_{\alpha_1}$ is a solution of (3.2) then choosing $\alpha_1$ such that $\gamma + \alpha_1 \theta > \beta = \alpha$, it follows from (3.5) that $y \in H_\alpha$.

Since the choice of $T_1$ can be done uniformly, the local solutions can be patch together to obtain existence on $[0, T]$. ∎

**Theorem 3.2. (Uniqueness).** Let $0 < \beta \le \alpha \le 1$, the function $F: R_+ \times R^m \times R^d \to \mathcal{L}(R^d, R^m)$ is the same as in Lemma 3.1 a) and c) with $\theta$- Hölder continuously differentiable by $y$ components, $0 < \theta \le 1$. If $0 < \beta - \gamma < \beta_0 < \theta\beta, 0 < \gamma$, then for any $x \in G_{\beta_0,\gamma} \cap H_\alpha([0,T], R^d)$ there exists a unique solution from $H_\beta([0,T], R^m)$ of equation (3.1).

Proof. Let $y, \tilde{y} \in H_\beta([0,T], R^m)$. The existence of $T_{x,F}y, T_{x,F}\tilde{y}$ in (3.3) follows from Lemma 3.1 a). For $a, b \in [0,T]$, $y, \tilde{y} \in B_T(\beta)$ from (3.4), applying Lemma 3.1 c)

$$|\int_a^b F(u, y(u), x(u)) dx(u) - \int_a^b F(u, \tilde{y}(u), x(u)) dx(u)|$$

$$\le \sum_{i=1}^m \sum_{j=1}^d |\int_a^b (F_{ij}(u, y(u), x(u)) - F_{ij}(u, \tilde{y}(u), x(u))) dx_j(u)|$$

$$\le \sum_{i=1}^m \sum_{j=1}^d (((|D_y F_{ij}, C| T^{\beta(1-\theta)}|y - \tilde{y}, H_\beta| + (|D_y F_{ij}, H_\theta(y)|(|y, H_\beta| + |\tilde{y}, H_\beta|)$$

$$+ |D_y F_{ij}, H_\theta(t)|T^{\theta(1-\beta)} + |D_y F_{ij}, H_\theta(x)||x, H_\alpha|T^{\theta(\alpha-\beta)})|y - \tilde{y}, C|)|x, G_{\beta\theta,\gamma}(x, [0,T])| |b-a|^{\gamma+\beta\theta}$$

$$+ 8|D_y F_{ij}, C| |x, H_\alpha([0,T])| |b-a|^\alpha |y - \tilde{y}, C|)$$

$$\le C(|b-a|^{\gamma+\beta\theta} + |x, H_\alpha([0,T])| |b-a|^\alpha)|y - \tilde{y}, H_\beta|,$$

where $C$ does not depend on $y, \tilde{y} \in B_T(\beta)$ and bounded by $T \in (0,1]$, since $|y - \tilde{y}, C| \le T^\beta |y - \tilde{y}, H_\beta|$. Thus, for the map $T_{x,F}$

$$|T_{x,F}y - T_{x,F}\tilde{y}, H_\beta([0,T])| \leq C\big(T^{\gamma-\beta(1-\theta)} + |x, H_\alpha([0,T])| T^{\alpha-\beta}\big)|y - \tilde{y}, H_\beta|.$$

Repeating now the arguments in the end of Theorem 3.1, we obtain for small enough $T = T_1, \beta < \alpha$ and subsequently for $\beta = \alpha$

$$|T_{x,F}y - T_{x,F}\tilde{y}, H_\beta([0,T_1])| \leq \frac{1}{2}|y - \tilde{y}, H_\beta([0,T_1])|.$$

Therefore, $T_{x,F}$ admits a unique fixed point in $B_{T_1}(\beta)$ which is the unique solution to (3.1) on the interval $[0, T_1], T_1 \leq T$. Patching together local solutions we obtain the unique solution on [0, T]. ∎

**Theorem 3.3. (Continuity of the Ito-Lyons map).** Let $y, \tilde{y} \in H_\beta, \beta > 0$, be unique solutions of equation (3.1) with inputs $(x, F, y_0), (\tilde{x}, \tilde{F}, \tilde{y}_0)$ such that $F(t, y, x), \tilde{F}(t, y, x)$ are the same as in Lemma 3.1 a) and c) with $\theta$- Hölder continuously differentiable by $t, y, x$ components, $0 < \theta \leq 1$, and $x, \tilde{x} \in G_{\beta_0, \gamma} \cap H_\alpha([0,T], \mathrm{R}^d), 0 < \beta - \gamma < \beta_0 < \theta\beta, 0 < \gamma$. Then
$|y - \tilde{y}, H_\beta([0,T])|$
$$\leq C(|y_0 - \tilde{y}_0| + |\tilde{x} - x, G_{\beta,\gamma}([0,T])| + |x - \tilde{x}, H_\alpha| + |x - \tilde{x}, \mathcal{C}|$$
$$+ |F - \tilde{F}, H_\beta([0,T])| + |F - \tilde{F}, \mathcal{C}|), \qquad (3.6)$$

where $C$ depends only on inputs.
Proof. If $y, \tilde{y}$ are solutions of (3.1) then from Theorem 3.1
$$|y, \mathcal{C}([0,T])| + |y, H_\beta([0,T])| \leq M,$$
$$|\tilde{y}, \mathcal{C}([0,T])| + |\tilde{y}, H_\beta([0,T])| \leq M < \infty,$$
where the constant $M$ depends only on equations inputs. Furthermore, applying estimates from Theorems 3.1 and 3.2, we obtain $(F, \tilde{F} \in H_1)$
$|y(a) - \tilde{y}(a) - (y(b) - \tilde{y}(b))|$
$$= |\int_a^b F(u, y(u), x(u))dx(u) - \int_a^b \tilde{F}(u, \tilde{y}(u), \tilde{x}(u))d\tilde{x}(u)|$$
$$\leq |\int_a^b \tilde{F}(u, \tilde{y}(u), \tilde{x}(u))d(\tilde{x}(u) - x(u)) + \int_a^b \big(\tilde{F}(u, \tilde{y}(u), \tilde{x}(u)) - \tilde{F}(u, y(u), \tilde{x}(u))\big)dx(u)$$
$$+ \int_a^b \big(\tilde{F}(u, y(u), \tilde{x}(u)) - \tilde{F}(u, y(u), x(u))\big)dx(u)$$
$$+ \int_a^b (\tilde{F}(u, y(u), x(u)) - F(u, y(u), x(u)))dx(u)|$$
$$\leq C(|\tilde{x} - x, G_{\beta,\gamma}([0,T])| |b-a|^{\gamma+\beta} + |\tilde{x} - x, H_\alpha||b-a|^\alpha$$
$$+ (|y - \tilde{y}, H_\beta([0,T])| + |y - \tilde{y}, \mathcal{C}|)|b-a|^{\gamma+\beta\theta} + |y - \tilde{y}, \mathcal{C}||b-a|^\alpha$$
$$+ (|x - \tilde{x}, H_\alpha([0,T])| + |x - \tilde{x}, \mathcal{C}|)|b-a|^{\gamma+\beta\theta} + |x - \tilde{x}, \mathcal{C}| |b-a|^\alpha$$
$$+ |F - \tilde{F}, H_\beta([0,T])| |b-a|^{\gamma+\beta} + |F - \tilde{F}, \mathcal{C}| |b-a|^\alpha),$$
where the constant $C$ depends only on inputs. Thus,
$$|y - \tilde{y}, H_\beta([0,T])| \leq C(|y_0 - \tilde{y}_0| + |y - \tilde{y}, H_\beta([0,T])|(T^{\gamma-\beta(1-\theta)} + T^{\alpha-\beta})$$
$$+ |F - \tilde{F}, H_\beta([0,T])| + |F - \tilde{F}, \mathcal{C}| + |x - \tilde{x}, H_\alpha| + |x - \tilde{x}, \mathcal{C}| + |\tilde{x} - x, G_{\beta,\gamma}([0,T])|).$$
Choosing $T = T_1$ such that $C(T_1^{\gamma-\beta(1-\theta)} + T_1^{\alpha-\beta}) \leq 1/2$ we obtain (3.6) for small enough $T_1$ and $\alpha > \beta$. Then (3.6) can be proved for $\alpha = \beta$ and, by using local solutions, for all $T$.

## 4. Remaining proofs

1. Integral existence. Proof of Theorem 2.1.b. Let $m > r > 0$ and $0 \leq a < b \leq 1$. From Definition 2.1, if $a_r = 2^{-r}n_r(a), b_r = 2^{-r}(n_r(b) + 1)$, then for $k \geq r, n_{k+1}(a_r) = 2n_k(a_r), n_{k+1}(b_r) = 2n_k(b_r) + 1$, and
$$\int_{\gamma_m(a,b,g)} f(t,x)dx - \int_{\gamma_r(a,b,g)} f(t,x)dx$$

$$= \sum_{k=r}^{m-1} \left( \int_{\gamma_{k+1}(a_r,b_r,g)} f(t,x)dx - \int_{\gamma_k(a_r,b_r,g)} f(t,x)dx \right) + \int_{\gamma_m(a,b,g)\setminus \gamma_m(a_r,b_r,g)} f(t,x)dx,$$

$$\int_{\gamma_{k+1}(a_r,b_r,g)} f(t,x)dx - \int_{\gamma_k(a_r,b_r,g)} f(t,x)dx$$

$$= \sum_{n=n_k(a_r)+1}^{n_k(b_r)+1} \left( \int_{h_{k,n-1}}^{h_{k+1,2(n-1)}} (f(t_{k,n-1},x) - f(t_{k+1,2(n-1)+1},x))dx \right.$$

$$\left. - \int_{h_{k+1,2(n-1)+1}}^{h_{k,n-1}} (f(t_{k+1,2(n-1)+1},x) - f(t_{k,n},x))dx \right)$$

$$- \int_{h_{k,n_k(a_r)}}^{h_{k+1,n_{k+1}(a_r)}} f(a_r,x)dx - \int_{h_{k+1,n_{k+1}(b_r)}}^{h_{k,n_k(b_r)}} f(b_r,x)dx, \qquad (4.1)$$

since

$$\sum_{n=n_{k+1}(a_r)+1}^{n_{k+1}(b_r)} \int_{h_{k+1,n-1}}^{h_{k+1,n}} f(t_{k+1,n},x)dx$$

$$= \sum_{n=n_k(a_r)+1}^{n_k(b_r)+1} \int_{h_{k+1,2(n-1)}}^{h_{k+1,2(n-1)+1}} f(t_{k+1,2(n-1)+1},x)dx + \sum_{n=n_k(a_r)+1}^{n_k(b_r)} \int_{h_{k+1,2(n-1)+1}}^{h_{k+1,2n}} f(t_{k+1,2n},x)dx,$$

$$\int_{h_{k+1,2n-1}}^{h_{k+1,2n}} f(t_{k+1,2n},x)dx - \int_{h_{k,n-1}}^{h_{k,n}} f(t_{k,n},x)dx$$

$$= \int_{h_{k+1,2n-1}}^{h_{k,n-1}} f(t_{k+1,2n},x)dx + \int_{h_{k,n}}^{h_{k+1,2n}} f(t_{k+1,2n},x)dx.$$

For $f(t,x) = f(t)$

$$\int_{\gamma_{k+1}(a_r,b_r,g)} f(t)dx - \int_{\gamma_k(a_r,b_r,g)} f(t)dx$$

$$= \sum_{n=n_k(a_r)+1}^{n_k(b_r)+1} \left( f(t_{k,n-1}) - 2f(t_{k+1,2(n-1)+1}) + f(t_{k,n}) \right) \left( h_{k+1,2(n-1)} - h_{k,n-1} \right)$$

$$- f(a_r)\left(h_{k+1,n_{k+1}(a_r)} - h_{k,n_k(a_r)}\right) + f(b_r)\left(h_{k+1,n_{k+1}(b_r)} - h_{k,n_k(b_r)}\right), \qquad (4.2)$$

since

$$h_{k+1,2(n-1)} - h_{k,n-1} = h_{k,n-1} - h_{k+1,2(n-1)+1}.$$

We omitted argument $g$ for all $h_{k,n}$ from Definition 2.1.

The next estimate provides the base for necessary and sufficient conditions for the existence of the introduced integral in Theorem 2.1.b.

**Estimate 4.1.** Let $0 \le a < b \le 1, m > r > k_0, 2^{-k_0-1} \le (b-a) \le 2^{-k_0+1}, |f, CH_\beta(K_g)| < \infty, \beta > 0$. Then integral (2.1) exists for any $g \in H_\alpha, \alpha > 0$, if

$$\mu(g,a,b,\beta) = \sum_{k=k_0+1}^{\infty} 2^{-(k+1)\beta+1} \sum_{n, [t_{k,n-1},t_{k,n}] \subset [a,b]} |h_{k+1,2(n-1)}(g) - h_{k+1,2(n-1)+1}(g)| < \infty,$$

and

$$\left| \int_a^b f(t,g(t))dg(t) - \int_{\gamma_{k_0+1}(a,b,g)} f(t,x)dx \right| \le |f, CH_\beta(K_g)|\mu(g,a,b,\beta) + C|f, C(K_g)| \, |g, H_\alpha| \, |b-a|^\alpha,$$

where $C$ depends only on $\alpha$.

Proof. Let $2^{-k}n_k(a) = a_k, 2^{-k}(n_k(b)+1) = b_k, k \ge r, n_k(a), n_k(b)$ are from Definition 2.1. With the presentation (4.1)

$$\left| \int_{\gamma_m(a,b,g)} f(t,x)dx - \int_{\gamma_r(a,b,g)} f(t,x)dx \right|$$

$$\leq \sum_{k=r}^{m-1} \sum_{n=n_k(a_r)+1}^{n_k(b_r)+1} \left( \int_{h_{k,n-1}(g)}^{h_{k+1,2(n-1)}(g)} |f(t_{k,n-1},x) - f(t_{k+1,2(n-1)+1},x)|dx \right.$$

$$\left. + \int_{h_{k+1,2(n-1)+1}(g)}^{h_{k,n-1}(g)} |f(t_{k+1,2(n-1)+1},x) - f(t_{k,n},x)|dx \right)$$

$$+ \left| \int_{h_{r,n_r(a_r)}(g)}^{h_{m,n_m(a_r)}(g)} f(a_r,x)dx \right| + \left| \int_{h_{m,n_m(b_r)}(g)}^{h_{r,n_r(b_r)}(g)} f(b_r,x)dx \right| + \left| \int_{\gamma_m(a,b,g)\backslash\gamma_m(a_r,b_r,g)} f(t,x)dx \right|$$

$$\leq |f, CH_\beta(K_g)| \sum_{k=r}^{m} 2^{-(k+1)\beta+1} \sum_{n,[t_{k,n-1},t_{k,n}]\subset[a,b]} |h_{k,n-1}(g) - h_{k+1,2(n-1)+1}(g)|$$

$$+ 4|f, C(K_g)||g, H_\alpha| 2^{-\alpha r} + \left| \int_{\gamma_m(a,b,g)\backslash\gamma_m(a_r,b_r,g)} f(t,x)dx \right|,$$

since

$$\left| \int_a^b (f(t,x) - f(t+h,x))dx \right| \leq |f, CH_\beta(K_g)||h|^\beta|a-b|,$$

and

$$|h_{m,n_m(a_r)}(g) - h_{r,n_r(a_r)}(g)| + |h_{m,n_m(b_r)}(g) - h_{r,n_r(b_r)}(g)| \leq 4|g, H_\alpha| 2^{-\alpha r}.$$

To estimate

$$B = \int_{\gamma_m(a,b,g)\backslash\gamma_m(a_r,b_r,g)} f(t,x)dx,$$

let us first consider the point $a$ in $\gamma_m(a,b,g)$. For $\gamma_m(a_m,a_r,g)$, we need to construct a special sequence of intervals $\{[a_{r_s+1}, a_{r_s}]\}_{s=0}^{s_1-1}$ since, in general, not all intervals $J_{k,n} = [n2^{-k}, (n+1)2^{-k}]$ are included in $[a_m, a_r]$ together with their parents $J_{k-1,[n/2]}$.

If $a_r = a_m$ then $B = 0$.

For $a_r > a_m$ define $r_0 = \max\{k \geq r; a_r = a_k\}$. According to Definition 2.1, for $r_0 \geq k \geq r$, there would be no extension of the interval $[a_r, b_r]$ to the left, and the first one would be the interval $[a_{r_0+1}, a_{r_0+1} + 2^{-r_0}]$ where $a_{r_0+1} + 2^{-r_0} = a_{r_0} = a_r$.

In the case of $a_{r_0+1} = a_m$, we define $[a_{r_0+1}, a_{r_0}]$ and stop further processing.

Otherwise, $a_{r_0+1} > a_m$. At this step, we already have the first element $[a_{r_0+1}, a_{r_0}]$ to construct the sequence, and we should repeat the previous arguments to create the next segment. Let's assume that on step $s_1 \geq 1$ we already created segments $\{[a_{r_s+1}, a_{r_s}]\}_{s=0}^{s_1-1}$. If $a_{r_{s_1-1}+1} = a_m$, we stop processing and use these intervals for further evaluations. Otherwise, $a_{r_{s_1-1}+1} > a_m$, and we can extend $a_{r_{s_1-1}+1}$ to the left by intervals with size less than $2^{-r_{s_1-1}}$. Define $r_{s_1} = \max\{k \geq r_{s_1-1}+1; a_{r_{s_1-1}+1} = a_k\}$, which gives us the next interval $[a_{r_{s_1}+1}, a_{r_{s_1}+1} + 2^{-r_{s_1}}]$, where $a_{r_{s_1}+1} + 2^{-r_{s_1}} = a_{r_{s_1-1}}$. Thus, in the finite number of steps $s_1 \geq 1$, we obtain the required sequence of segments

$$\{[a_{r_s+1}, a_{r_s}]\}_{s=0}^{s_1-1}, r_{s_1} = m \geq r_{s_1-1}+1 > r_{s_1-1} > .. > r_1 \geq r_0+1 > r_0 \geq r,$$

with segments length $2^{-r_s}$ and $a_{r_{s_1-1}+1} = a_m$.

With these notations and arguments similar to those used in the previous estimates,

$$\left| \int_{\gamma_m(a_m,a_{r_0},g)} f(t,x)dx \right|$$

$$= \left| \sum_{q=1}^{s_1} \left( \int_{\gamma_m(a_{r_q},a_{r_{q-1}},g)} - \int_{\gamma_{r_{q-1}+1}(a_{r_q},a_{r_{q-1}},g)} \right) f(t,x)dx \right|$$

$$+\sum_{q=1}^{s_1}\int_{\gamma_{r_{q-1}+1}(a_{r_q},a_{r_{q-1}},g)}f(t,x)dx+\sum_{q=1}^{s_1}\int_{h_{m,2^{m-r_{q-1}}n_{r_{q-1}}(a)-1}}^{h_{m,2^{m-r_{q-1}}n_{r_{q-1}}(a)}}f\left(a_{r_{q-1}},x\right)dx|$$

$$\leq |f,\mathcal{CH}_\beta(K_g)|\sum_{q=1}^{s_1}\sum_{k=r_{q-1}+1}^{m}2^{-(k+1)\beta+1}\sum_{n,[t_{k,n-1},t_{k,n}]\subset[a_{r_q},a_{r_{q-1}}]}|h_{k+1,2(n-1)}-h_{k+1,2(n-1)+1}|$$

$$+8|f,\mathcal{C}(K_g)||g,H_\alpha|\sum_{q=1}^{s_1}(2^{-\alpha(r_{q-1}+1)}+2^{-\alpha m})$$

$$\leq |f,H_\beta(K_g)|\sum_{k=r_0+1}^{m}2^{-(k+1)\beta+1}\sum_{n,[t_{k,n-1},t_{k,n}]\subset[a_m,a_{r_0}]}|h_{k+1,2(n-1)}-h_{k+1,2(n-1)+1}|$$

$$+16|f,\mathcal{C}(K_g)||g,H_\alpha|\sum_{q=1}^{s_1}2^{-\alpha(r_0+q)},$$

since
$$r_{s_1}=m;\ a_{r_q}=a_{r_{q-1}+1}=a_{r_{q-1}}-2^{-r_{q-1}},$$
$$2^{-m}2^{m-r_{q-1}}n_{r_{q-1}}(a)=a_{r_{q-1}},\ 2^{-m}(2^{m-r_{q-1}}n_{r_{q-1}}(a)-1)=a_{r_{q-1}}-2^{-m}.$$

With the same estimate for point $b$ in $\gamma_m(a,b,g)$ and two previous estimates, we obtain

$$|\int_{\gamma_m(a,b,g)}f(t,x)dx-\int_{\gamma_r(a,b,g)}f(t,x)dx|\to 0, when\ m,r\to\infty,$$

which means that

$$\{\int_{\gamma_m(a,b,g)}\omega\}$$

is a Cauchy sequence; thus, it converges to a number that is the value of integral (2.1).
To prove (2.3), we should choose $r=k_0$ and use estimates from Theorem 2.1. ∎

To proof Theorem 2.8 we need to repeat the arguments of Estimate 4.1 but use the presentation (4.2) instead of (4.1).
For $f\in\mathbb{G}_{g,f'}$, by Definition 2.3
$$|f(t_{k,n-1})-2f(t_{k+1,2(n-1)+1})+f(t_{k,n})|\leq C\ \psi(2^{-k})|f',H_{\psi(\delta)}([0,1])||h_{k+1,2n}(g)-h_{k+1,2n+1}(g)|,$$
where $C$ does not depend on $f$ and $f'$. From presentation (4.2)

$$|\sum_{n=n_k(a_r)+1}^{n_k(b_r)}\left(f(t_{k,n-1})-2f(t_{k+1,2(n-1)+1})+f(t_{k,n})\right)(h_{k+1,2(n-1)}-h_{k,n-1})|$$

$$\leq C\ |f',H_{\psi(\delta)}([0,1])|\psi(2^{-k})\sum_{n=n_k(a_r)+1}^{n_k(b_r)}|h_{k+1,2n}(g)-h_{k+1,2n+1}(g)|^2. ∎$$

In Theorem 2.1, we showed that condition (2.2) is sufficient for the existence of integral (2.1) and therefore for (2.4).
In the next theorem, we prove that it is also necessary for the existence of integral (2.4).

**Theorem 4.1.** (Theorem 2.2 necessity). Let $g\in\mathcal{C}, 1>\beta>0$, and integral (2.4) with respect to the function $g$ exists for any function from class $H_\beta$; then, the function $g$ satisfies condition (2.2).

In the proof, we construct a function $f\in H_\beta$ such that

$$\infty > \int_0^1 f(t)dg(t)=\sum_{k=1}^{\infty}2^{-\beta(k+1)}|h_{k+1,2n}(g)-h_{k+1,2n+1}(g)|,$$

which proves condition (2.2).
For each binary interval $I_{k,n}=[n2^{-k},(n+1)2^{-k}], k>0, 0\leq n<2^k$, define the continuous function

$f_{k,n}:[0,1] \to R$ as $f_{k,n}|_{[0,1]\setminus I_{k,n}} = 0$, $f_{k,n}(t_{k+1,2n+1}) = 2^{-(k+1)\beta}\,sign(h_{k+1,2n+1} - h_{k+1,2n})$, and $f_{k,n}$ is the linear extension on [0,1]. Finally, define $f:[0,1] \to R$ as

$$f(t) = \sum_{k=1}^{\infty}\sum_{n=0}^{2^k-1} f_{k,n}(t).$$

It follows from Lemma 5.5 of Guseynov (2016) that $f \in H_\beta([0,1])$, and by the theorem's assumption,

$$\lim_{r\to\infty}\int_{\gamma_r(0,1,g)} f(t)dg(t) = \int_0^1 f(t)dg(t)$$

exists and is finite.

Let $k_0 > 1$. By the construction of functions $f_{k,m}$, $f_{k,m}(t_{k_0,n}) = 0$ for $k \geq k_0$, and for $1 \leq k < k_0$

$$f_{k,m}(t_{k_0,i}) = \begin{cases} 0, & i2^{-k_0} \leq m2^{-k}; \\ a_{k,m}2^{k+1}(i2^{-k_0} - m2^{-k}), & m2^{-k} < i2^{-k_0} \leq m2^{-k} + 2^{-k-1}; \\ a_{k,m}2^{k+1}((m+1)2^{-k} - i2^{-k_0}), & m2^{-k} + 2^{-k-1} < i2^{-k_0} < (m+1)2^{-k}; \\ 0, & i2^{-k_0} \geq (m+1)2^{-k}, \end{cases}$$

where $a_{k,m} = 2^{-(k+1)\beta}\,sign(h_{k+1,2m+1} - h_{k+1,2m})$. Having this and Definition 2.1,

$$\int_{\gamma_{k_0}(0,1,g)} f(t)dx = \sum_{n=1}^{2^{k_0}-1} f(t_{k_0,n})(h_{k_0,n} - h_{k_0,n-1})$$

$$= \sum_{k=1}^{k_0-1}\sum_{m=0}^{2^k-1}\sum_{i=m2^{k_0-k}+1}^{(m+1)2^{k_0-k}-1} f_{k,m}(t_{k_0,i})(h_{k_0,i} - h_{k_0,i-1}),$$

and

$$\sum_{i=m2^{k_0-k}+1}^{(m+1)2^{k_0-k}-1} f_{k,m}(t_{k_0,i})(h_{k_0,i} - h_{k_0,i-1})$$

$$= a_{k,m}\Big(\sum_{i=m2^{k_0-k}+1}^{m2^{k_0-k}+2^{k_0-k-1}-1} 2^{-(k_0-k-1)}(i - m2^{k_0-k})(h_{k_0,i} - h_{k_0,i-1})$$
$$+ (h_{k_0,2^{k_0-k-1}(2m+1)} - h_{k_0,2^{k_0-k-1}(2m+1)-1})$$
$$+ \sum_{i=m2^{k_0-k}+2^{k_0-k-1}+1}^{(m+1)2^{k_0-k}-1} 2^{-(k_0-k-1)}((m+1)2^{k_0-k} - i)(h_{k_0,i} - h_{k_0,i-1})\Big)$$

$$= a_{k,m}\Big(\sum_{i=1}^{2^{k_0-k-1}-1} 2^{-(k_0-k-1)}i * \Big(h_{k_0,i+m2^{k_0-k}} - h_{k_0,i+m2^{k_0-k}-1}\Big) - h_{k_0,2^{k_0-k-1}(2m+1)-1}$$

$$+ h_{k_0,2^{k_0-k-1}(2m+1)} + \sum_{i=1}^{2^{k_0-k-1}-1} 2^{-(k_0-k-1)}i\Big(h_{k_0,(m+1)2^{k_0-k}-i} - h_{k_0,(m+1)2^{k_0-k}-i-1}\Big)\Big)$$

$$= a_{k,m}\Big(-\sum_{i=0}^{2^{k_0-k-1}-1} 2^{-(k_0-k-1)} * h_{k_0,i+m2^{k_0-k}} + \sum_{i=1}^{2^{k_0-k-1}} 2^{-(k_0-k-1)} h_{k_0,(m+1)2^{k_0-k}-i}\Big)$$

$$= a_{k,m}2^{k+1}\Big(-\sum_{i=0}^{2^{k_0-k-1}-1}\int_{m2^{-k}+i2^{-k_0}}^{m2^{-k}+(i+1)2^{-k_0}} g(t)dt + \sum_{i=1}^{2^{k_0-k-1}}\int_{(m+1)2^{-k}-i2^{-k_0}}^{(m+1)2^{-k}-(i-1)2^{-k_0}} g(t)dt\Big)$$

$$= 2^{-(k+1)\beta}|h_{k+1,2m+1} - h_{k+1,2m}|.$$

Thus,

$$\int_{\gamma_{k_0}(0,1,g)} f(t)dg(t) = \sum_{k=1}^{k_0-1}\sum_{m=0}^{2^k-1} 2^{-(k+1)\beta}|h_{k+1,2m+1} - h_{k+1,2m}|. \blacksquare$$

2. **Integration in $\alpha$-Hölder classes with $\alpha \leq 2^{-1}$.** Denote by $\mathcal{C}$ the space of all real-valued continuous functions $x:[0,1] \to R$. We shall use the $\sigma$-additive measure W on $\mathcal{C}$ introduced by Wiener (1930), which is defined through the Gaussian measure on the cylinder sets

$$C(t_1, \ldots, t_n; I_1, \ldots, I_n) = \{x \in \mathcal{C}; x(t_k) \in I_k, 1 \leq k \leq n\},$$

$$W\big(C(t_1, \ldots, t_n; I_1, \ldots, I_n)\big) = \int_{I_1}\int_{I_2}\cdots\int_{I_n} \prod_{k=1}^{n} \frac{d\xi_k}{\sqrt{2\pi(t_k - t_{k-1})}} \exp\left(-\frac{(\xi_k - \xi_{k-1})^2}{2(t_k - t_{k-1})}\right),$$

where $0 < t_1 < \cdots < t_n \leq 1$, $I_k = (a_k, b_k)$, $1 \leq k \leq n$ are real intervals and $t_0 = \xi_0 = 0$.

Let $(\mathcal{C}, \mathcal{F}, W)$ be a probability (Wiener) space, where $\mathcal{F}$ is the $\sigma$-algebra obtained from the cylinders by the operations of countable unions, intersections, and the operation of complement, and let

$$\int_\mathcal{C} \mathcal{L}(x) dW(x)$$

be a Lebesgue integral defined on the space $(\mathcal{C}, \mathcal{F}, W)$ for any measurable function $\mathcal{L}: \mathcal{C} \to R$. The finite dimensional example of the functional $\mathcal{L}$ is a continuous function $F: R^n \to R$ that depends only on a finite number of values of $x$, i.e.,

$$\mathcal{L}(x) = F\big(x(t_1), \ldots, x(t_n)\big), 0 < t_1 < \cdots < t_n \leq 1,$$

and by Wiener (1930), the corresponding integral can be expressed as an n-fold Riemann integral

$$\int_\mathcal{C} \mathcal{L}(x) dW(x)$$

$$= \int_{-\infty}^{\infty}\int_{-\infty}^{\infty}\cdots\int_{-\infty}^{\infty} F(\xi_1, \ldots, \xi_n) \prod_{k=1}^{n} \frac{d\xi_k}{\sqrt{2\pi(t_k - t_{k-1})}} \exp\left(-\frac{(\xi_k - \xi_{k-1})^2}{2(t_k - t_{k-1})}\right), \xi_0 = t_0 = 0.$$

We use the Wiener measure to evaluate the deterministic integral over Wiener process sample paths or BMSPs. We also use the

**Borel-Cantelli Lemma.** Let $\{A_n\}$ be a sequence of sets in a probability space $(\Omega, \Sigma, P)$. Then, the set

$$\{A_n \text{ infinitely often }\} = \overline{\lim_{n \to \infty}} A_n = \bigcap_{n=1}^{\infty}\bigcup_{k=n}^{\infty} A_k$$

consists of the elements $x \in \Omega$ such that x belongs to infinitely many sets from the sequence $\{A_n\}$, and if

$$\sum_{n=1}^{\infty} P(A_n) < \infty, \text{ then } P(\{A_n \ i.o.\}) = 0.$$

For further evaluations, we need to approximate $h_{k+1,2n+1} - h_{k+1,2n}$ through independent increments of function $g$.

**Lemma 4.1.** Approximation $h_{k+1,2n} - h_{k+1,2n+1}$. Let $l, k > 0$, $q = l + k + 1$, $n = 0, \ldots, 2^k - 1$ and $g$ be continuous on $[0,1]$; then,

$$(Approx)(h_{k+1,2n} - h_{k+1,2n+1})$$

$$= (Approx)(2^{k+1}\int_{t_{k+1,2n}}^{t_{k+1,2n+1}} g(\tau)d\tau - 2^{k+1}\int_{k+1,2n+1}^{t_{k+1,2(n+1)}} g(\tau)d\tau)$$

$$= ((l)Approx)(h_{k+1,2n} - h_{k+1,2n+1}) = \Delta^\theta_{k+1,2n+1,l}g + 2^{-l}\Delta^2_{k+1,2n}g,$$

where

$$\Delta^\theta_{k+1,2n+1,l}g = \sum_{p=2n2^l}^{2(n+1)2^l-1} (g_{q,p} - g_{q,p+1})\theta(2^{-l}p), g_{q,p} = g(t_{q,p}), t_{q,p} = p * 2^{-q},$$

$$\Delta^2_{k+1,2n}g = g_{k+1,2n} - 2g_{k+1,2n+1} + g_{k+1,2(n+1)},$$

$\theta(t) = frac(t), 2n \leq t < 2n + 1, \theta(t) = 1 - frac(t), 2n + 1 \leq t < 2(n + 1),$

$$((l)Approx)(h_{k+1,2n} - h_{k+1,2n+1}) \to (h_{k+1,2n} - h_{k+1,2n+1}), l \to \infty,$$

for all $g \in \mathcal{C}$, and

$$\int_{\mathcal{C}} \Delta^{\theta}_{k+1,2n+1,l} x \, \Delta^{\theta}_{k+1,2m+1,l} x \, dW(x) = 0, n \neq m,$$

$$\int_{\mathcal{C}} (\Delta^{\theta}_{k+1,2n+1,l} x)^2 dW(x) = 2^{-(k+1)} \sum_{p=2n2^l}^{2(n+1)2^l-1} 2^{-l}(\theta(2^{-l}p))^2 = 2^{-k-1} C_l(\theta),$$

$$C_l(\theta) \to \frac{2}{3}, l \to \infty.$$

Proof. For $t_0 = \tau_0, t_1 = \tau_s, \tau_{n+1} - \tau_n = (t_1 - t_0)/s$, we can represent the integral sum for $g \in \mathcal{C}$ as

$$(Approx)(\frac{1}{t_1 - t_0} \int_{t_0}^{t_1} g(\tau) d\tau) = \frac{1}{t_1 - t_0} \sum_{n=0}^{s-1} g(\tau_n)(\tau_{n+1} - \tau_n)$$

$$= \frac{1}{t_1 - t_0} \sum_{n=1}^{s} (g(\tau_{n-1}) - g(\tau_n))(\tau_n - \tau_s) + g(\tau_0),$$

and

$$((l)Approx)(h_{k+1,2n} - h_{k+1,2n+1})$$

$$= 2^{k+1} \sum_{p=0}^{2^l-1} \left(g_{q,2n2^l+p} - g_{q,2n2^l+p+1}\right) (2^{-q}(2n2^l + p + 1) - (2n+1)2^{-(k+1)}) + g_{k+1,2n}$$

$$- 2^{k+1} \sum_{p=0}^{2^l-1} \left(g_{q,(2n+1)2^l+p} - g_{q,(2n+1)2^l+p+1}\right) (2^{-q}((2n+1)2^l + p + 1))$$

$$- 2(n+1)2^{-(k+1)}) - g_{k+1,2n+1}$$

$$= \sum_{p=2n2^l}^{(2n+1)2^l-1} (g_{q,p} - g_{q,p+1}) (2^{-l}(p+1) - (2n+1)) + g_{k+1,2n}$$

$$- \sum_{p=(2n+1)2^l}^{2(n+1)2^l-1} (g_{q,p} - g_{q,p+1}) (2^{-l}(p+1) - 2(n+1)) - g_{k+1,2n+1}$$

$$= \sum_{p=2n2^l}^{(2n+1)2^l-1} (g_{q,p} - g_{q,p+1}) (frac(2^{-l}p) + 2^{-l} - 1) + g_{k+1,2n}$$

$$- \sum_{p=(2n+1)2^l}^{2(n+1)2^l-1} (g_{q,p} - g_{q,p+1}) (frac(2^{-l}p) - 1 + 2^{-l}) - g_{k+1,2n+1}$$

$$= \sum_{p=2n2^l}^{2(n+1)2^l-1} (g_{q,p} - g_{q,p+1}) \theta(2^{-l}p) + 2^{-l}(g_{k+1,2n} - 2g_{k+1,2n+1} + g_{k+1,2(n+1)}).$$

Since function $g \in \mathcal{C}$,

$$\frac{1}{|t_1 - t_0|} | \int_{t_0}^{t_1} g(\tau) d\tau - \sum_{n=0}^{s-1} g(\tau_n)(\tau_{n+1} - \tau_n)| \leq \sup_{|t-\tau|<|t_1-t_0|/s} |g(t) - g(\tau)| \to 0, s \to \infty. \blacksquare$$

**Proof of Theorem 2.7**. Let $a = (\frac{2}{3\pi})^{\frac{1}{2}}$ and $k > 1$. We shall evaluate the integral

$$\int_{\mathcal{C}} (2^{-k/2} \sum_{m=0}^{2^k-1} |h_{k+1,2m+1}(x) - h_{k+1,2m}(x)| - a)^2 dW(x)$$

$$= \int_C 2^{-k} \sum_{m=0}^{2^k-1} (h_{k+1,2m+1}(x) - h_{k+1,2m}(x))^2 dW(x)$$

$$+ (2^{-\frac{k}{2}} \sum_{m=0}^{2^k-1} \int_C dW(x) |h_{k+1,2m+1}(x) - h_{k+1,2m}(x)| - a)^2$$

$$- 2^{-k} \sum_{m=0}^{2^k-1} (\int_C dW(x) |h_{k+1,2m+1}(x) - h_{k+1,2m}(x)|)^2. \quad (4.3)$$

To prove (4.3) we need to show that

$$\int_C |h_{k+1,2n+1}(x) - h_{k+1,2n}(x)| |h_{k+1,2m+1}(x) - h_{k+1,2m}(x)| dW(x)$$

$$= \int_C |h_{k+1,2n+1}(x) - h_{k+1,2n}(x)| dW(x) \int_C |h_{k+1,2m+1}(x) - h_{k+1,2m}(x)| dW(x), n \neq m.$$

Based on Lemma 4.1, we have it for $(Approx)(h_{k+1,2n} - h_{k+1,2n+1})$; then, it follows from the Fatou theorem about the convergence of a sequence of nonnegative functions.

**Lemma 4.2.** Let $k > 1$, $n = 0, \ldots, 2^k - 1$; then,

$$\int_C |h_{k+1,2n+1}(x) - h_{k+1,2n}(x)| dW(x) = 2^{-\frac{k}{2}} (\frac{2}{3\pi})^{\frac{1}{2}}.$$

Proof. We shall use Lemma 4.1 to approximate $h_{k+1,2n+1}(g) - h_{k+1,2n}(g)$. Let $l > 0$, $q = l + k + 1$. Consider

$$(Approx) \int_C |h_{k+1,2n+1}(x) - h_{k+1,2n}(x)| dW(x)$$

$$= \int_C | \sum_{p=2n2^l}^{2(n+1)2^l-1} (x_{q,p} - x_{q,p+1}) \theta(2^{-l}p) + 2^{-l} \Delta^2_{k+1,2n} x | dW(x)$$

$$= \int_{-\infty}^{\infty} \int_{-\infty}^{\infty} \cdots \int_{-\infty}^{\infty} | \sum_{r=2n2^l}^{2(n+1)2^l-1} (\xi_r - \xi_{r+1}) \theta(2^{-l}r) + 2^{-l} (\sum_{r=2n2^l}^{(2n+1)2^l-1} (\xi_r - \xi_{r+1}) - \sum_{r=(2n+1)2^l}^{2(n+1)2^l-1} (\xi_r - \xi_{r+1})) |$$

$$* \prod_{p=2n2^l-1}^{2(n+1)2^l-1} \frac{d\xi_{p+1}}{\sqrt{2\pi(t_{p+1} - t_p)}} \exp\left(-\frac{(\xi_p - \xi_{p+1})^2}{2(t_{p+1} - t_p)}\right), \xi_{2n2^l-1} = t_{2n2^l-1} = 0, t_p = p2^{-q}.$$

Let us first evaluate ($L = 2^{l+1} - 1$, $\theta(2^{-l}r) = 0|_{r=0}$)

$$B_1 = \int_{R_{L+1}} |\sum_{r=1}^{L} (\xi_r - \xi_{r+1}) \theta(2^{-l}r)| \prod_{p=0}^{L} \frac{d\xi_{p+1}}{\sqrt{2\pi(t_{p+1} - t_p)}} \exp\left(-\frac{(\xi_p - \xi_{p+1})^2}{2(t_{p+1} - t_p)}\right)$$

$$= \int_{R_L} \prod_{p=0}^{L-1} \frac{d\xi_{p+1}}{\sqrt{2\pi(t_{p+1} - t_p)}} \exp\left(-\frac{(\xi_p - \xi_{p+1})^2}{2(t_{p+1} - t_p)}\right)$$

$$* \int_{-\infty}^{\infty} |\sum_{r=1}^{L-1} (\xi_r - \xi_{r+1}) \theta(2^{-l}r) + (\xi_L - \xi_{L+1}) \theta(2^{-l}L)| \frac{d\xi_{L+1}}{\sqrt{2\pi(t_{L+1} - t_L)}} \exp\left(-\frac{(\xi_L - \xi_{L+1})^2}{2(t_{L+1} - t_L)}\right)$$

$$= 2^{-q/2} \int_{R_L} |\sum_{r=1}^{L} \xi_r| \prod_{p=1}^{L} \frac{d\xi_p}{\theta(2^{-l}p)\sqrt{2\pi}} \exp\left(-\frac{(\xi_p)^2}{2(\theta(2^{-l}p))^2}\right)$$

$$= 2 * 2^{-q/2} \int_{\sum_{r=1}^{L} \xi_r > 0} \sum_{r=1}^{L} \xi_r \prod_{p=1}^{L} \frac{d\xi_p}{\theta(2^{-l}p)\sqrt{2\pi}} \exp\left(-\frac{\xi_p^2}{2(\theta(2^{-l}p))^2}\right)$$

$$= 2 * 2^{-q/2} \sum_{r=1}^{L} \int_{R_{L-1}} \prod_{p=1, p \neq r}^{L} \left( \frac{d\xi_p}{\theta(2^{-l}p)\sqrt{2\pi}} \exp\left(-\frac{\xi_p^2}{2(\theta(2^{-l}p))^2}\right) \right)$$

$$* \int_{\xi_r > -\sum_{p=1, p \neq r}^{L} \xi_p} \frac{\xi_r d\xi_r}{\theta(2^{-l}r)\sqrt{2\pi}} \exp\left(-\frac{\xi_r^2}{2(\theta(2^{-l}r))^2}\right)$$

$$= \frac{2 * 2^{-q/2}}{\sqrt{2\pi}} \sum_{r=1}^{L} \int_{R_{L-1}} \prod_{p=1, p \neq r}^{L} \left( \frac{d\xi_p}{\theta(2^{-l}p)\sqrt{2\pi}} \exp\left(-\frac{\xi_p^2}{2(\theta(2^{-l}p))^2}\right) \right)$$

$$* \theta(2^{-l}r) \exp\left(-\frac{1}{2}\left(\frac{1}{\theta(2^{-l}r)} \sum_{p=1, p \neq r}^{L} \xi_p\right)^2\right).$$

To further evaluate $B_1$, for each $r = 1, 2, \ldots L$ consider

$$F_m = \int_{-\infty}^{\infty} \frac{d\xi_{p_m}}{\theta(2^{-l}p_m)\sqrt{2\pi}} \exp\left(-\frac{(\xi_{p_m})^2}{2(\theta(2^{-l}p_m))^2}\right) \frac{1}{a_{m-1}} \exp\left(-\frac{1}{2}\frac{(\xi_{p_m} + A_m)^2}{(a_{m-1})^2}\right),$$

where $m = 1, 2, \ldots L, p_m \neq r$,

$$A_m = \sum_{i=1}^{L-(m+1)} \xi_{p_i}, 1 \le p_i < p_j \le L, \text{ if } i < j$$

(the set $\{p_1, \ldots, p_{L-(m+1)}\}$ does not include elements $p_m$ and $r$), and $a_0 = \theta(2^{-l}r)$. Since

$$\frac{\xi_{p_m}^2}{(\theta(2^{-l}p_m))^2} + \frac{(\xi_{p_m} + A_m)^2}{(a_{m-1})^2} = \frac{(A_m)^2}{(a_{m-1})^2 + (\theta(2^{-l}p_m))^2}$$

$$+ \left(\frac{1}{(\theta(2^{-l}p_m))^2} + \frac{1}{(a_{m-1})^2}\right) \left(\xi_{p_m} + \left(\frac{1}{(\theta(2^{-l}p_m))^2} + \frac{1}{(a_{m-1})^2}\right)^{-1} \frac{1}{(a_{m-1})^2} A_m\right)^2,$$

then

$$F_m = \exp\left(-\frac{1}{2} \frac{(A_m)^2}{((a_{m-1})^2 + (\theta(2^{-l}p_m))^2)}\right) \frac{1}{a_{m-1}} \left(\frac{1}{(\theta(2^{-l}p_m))^2} + \frac{1}{(a_{m-1})^2}\right)^{-1/2} \frac{1}{\theta(2^{-l}p_m)}$$

$$= \exp\left(-\frac{1}{2} \frac{(A_m)^2}{(a_m)^2}\right) \frac{1}{a_m}.$$

where $(a_m)^2 = (a_{m-1})^2 + (\theta(2^{-l}p_m))^2$. Now, subsequently extracting an element from the set $\{1, 2, \ldots, L\} \setminus \{r\}$ with no return and applying the induction, we obtain ($A_{L-2} = \xi_{p_{L-1}}$)

$$B_1 = \frac{2 * 2^{-q/2}}{\sqrt{2\pi}} \sum_{r=1}^{L} (\theta(2^{-l}r))^2 \int_{R_{L-2}} \left( \prod_{p=1, p \neq r, p \neq p_1}^{L} \frac{d\xi_p}{\theta(2^{-l}p)\sqrt{2\pi}} \exp\left(-\frac{\xi_p^2}{2(\theta(2^{-l}p))^2}\right) \right) F_1$$

$$= \frac{2 * 2^{-q/2}}{\sqrt{2\pi}} \sum_{r=1}^{L} (\theta(2^{-l}r))^2 \int_{R_{L-2}} \left( \prod_{p=1, p \neq r, p \neq p_1}^{L} \frac{d\xi_p}{\theta(2^{-l}p)\sqrt{2\pi}} \exp\left(-\frac{\xi_p^2}{2(\theta(2^{-l}p))^2}\right) \right) \exp\left(-\frac{1}{2}\frac{(A_1)^2}{(a_1)^2}\right) \frac{1}{a_1}$$

$$= \frac{2 * 2^{-q/2}}{\sqrt{2\pi}} \sum_{r=1}^{L} (\theta(2^{-l}r))^2 \int_{R_1} \left(\frac{d\xi_{p_{L-1}}}{\theta(2^{-l}p_{L-1})\sqrt{2\pi}} \exp\left(-\frac{(\xi_{p_{L-1}})^2}{2(\theta(2^{-l}p_{L-1}))^2}\right)\right) \exp\left(-\frac{1}{2}\frac{(\xi_{p_{L-1}})^2}{(a_{L-2})^2}\right) \frac{1}{a_{L-2}}$$

$$= \frac{2 * 2^{-q/2}}{\sqrt{2\pi}} \sum_{r=1}^{L} (\theta(2^{-l}r))^2 \frac{1}{a_{L-1}} = \frac{2 * 2^{-q/2}}{\sqrt{2\pi}} \sum_{r=1}^{L} (\theta(2^{-l}r))^2 \left(\sum_{r=1}^{L} (\theta(2^{-l}r))^2\right)^{-\frac{1}{2}}$$

$$= \frac{2 * 2^{-(k+1)/2}}{\sqrt{2\pi}} \left(\sum_{r=1}^{L} 2^{-l}(\theta(2^{-l}r))^2\right)^{1/2} \to \frac{2^{-k/2}}{\sqrt{\pi}} \left(2\int_0^1 t^2 dt\right)^{1/2}, l \to \infty.$$

Similarly,

$$\int_{-\infty}^{\infty} \int_{-\infty}^{\infty} \ldots \int_{-\infty}^{\infty} \left| 2^{-l} \left( \sum_{p=2n2^l}^{(2n+1)2^l-1} (\xi_p - \xi_{p+1}) - \sum_{p=(2n+1)2^l}^{2(n+1)2^l-1} (\xi_p - \xi_{p+1}) \right) \right|$$

$$* \prod_{p=2n2^{l-1}}^{2(n+1)2^l-1} \frac{d\xi_{p+1}}{\sqrt{2\pi(t_{p+1}-t_p)}} \exp\left(-\frac{(\xi_p-\xi_{p+1})^2}{2(t_{p+1}-t_p)}\right) \le 2^{-l}C \to 0, l \to \infty,$$

which, with Lemma 4.1 and the Fatou theorem about the convergence of a sequence of nonnegative functions, completes the proof of Lemma 4.2.

From Lemmas 4.1 and 4.2 and (4.3),

$$\int_{\mathcal{C}} (2^{-k/2} \sum_{m=0}^{2^k-1} |h_{k+1,2m+1}(x) - h_{k+1,2m}(x)| - (\frac{2}{3\pi})^{\frac{1}{2}})^2 dW(x) = \frac{2}{3} 2^{-(k+1)} - \frac{2}{3\pi} 2^{-k}. \quad (4.4)$$

An argument similar to that in Paley et al. (1933), see also Cameron and Martin (1947), yields the W-a.e. convergence in Theorem 2.7. Indeed, consider sets

$$E_k = \{x \in \mathcal{C}; |\sigma_k(x)| \ge 2^{-k/4}\},$$

where

$$\sigma_k(x) = |2^{-k/2} \sum_{m=0}^{2^k-1} |h_{k+1,2m+1}(x) - h_{k+1,2m}(x)| - (\frac{2}{3\pi})^{\frac{1}{2}}|.$$

Then, by (4.4)

$$W(E_k) \le C' * 2^{-k/2}, \text{ hence } \sum_{k>0} W(E_k) < \infty,$$

and by the Borel-Cantelli lemma,

$$W(\bigcap_{n=1}^{\infty} \bigcup_{k=n}^{\infty} E_k) = 0.$$

Thus, for

$$x \in \mathcal{C} \setminus (\bigcap_{n=1}^{\infty} \bigcup_{k=n}^{\infty} E_k)$$

(hence for W-a.a. $x$) there exists $k_1(x)$ such that for any $k > k_1(x)$, $0 \le \sigma_k(x) < 2^{-k/4} \to 0$, when $k \to \infty$. ∎

3. Integral calculus. To prove (2.6) for $0 \le a < c < b \le 1$ and $k > 0$, consider

$$\int_{\gamma_k(a,b,g)} f(t,x)dx = \sum_{n=n_k(a)+1}^{n_k(b)} \int_{h_{k,n-1}}^{h_{k,n}} f(t_{k,n},x)dx = \sum_{n=n_k(a)+1}^{n_k^{[a,c]}(c)} \int_{h_{k,n-1}}^{h_{k,n}} f(t_{k,n},x)dx$$

$$+ \sum_{n=n_k^{[c,b]}(c)+1}^{n_k(b)} \int_{h_{k,n-1}}^{h_{k,n}} f(t_{k,n},x)dx + \sum_{n=n_k^{[a,c]}(c)+1}^{n_k^{[c,b]}(c)} \int_{h_{k,n-1}}^{h_{k,n}} f(t_{k,n},x)dx, \quad (4.5)$$

where $n_k^{[a,c]}(c)$ is $n_k(c)$ for interval $[a,c]$ and $n_k^{[c,b]}(c)$ is $n_k(c)$ for interval $[c,b]$ from Definition 2.1. Additionally,

$$|\sum_{n=n_k^{[a,c]}(c)+1}^{n_k^{[c,b]}(c)} \int_{h_{k,n-1}}^{h_{k,n}} f(t_{k,n},x)dx| \le 6 * |f, \mathcal{C}(K_g)||g, H_\alpha| 2^{-\alpha k}.$$

Thus, with condition (2.2), letting $k \to \infty$ in (4.5), we obtain (2.6).

The Green formula for Jordan domains with irregular boundary.

**Lemma 4.3.** Let $\gamma$ be a contour in the plane $(t,x)$ given by a parametric presentation $\gamma(t) = (t, g(t))$, $t \in [a,b] \subset [0,1]$, and $\theta$ be a segment connecting points $\gamma(a)$ and $\gamma(b)$ such that $\gamma \cup \theta$ be a closed Jordan curve, which is the boundary of a bounded domain $G$. Let the function $f$ be the same as in Theorem 2.1.b, $f \in$ ACL (Absolutely Continuous on Lines), $\partial f/\partial t$ is summable in $K_g$. Then, for $g \in GH_\beta \cap H_\alpha, \alpha, \beta > 0$ or if $\beta > 2^{-1}$ for W-a.a. $g \in \mathcal{C}$

$$\int_{\partial G} f(t,x)dx = \int_G \frac{\partial}{\partial t} f(t,x)dtdx. \tag{4.6}$$

Proof. Let
$$\gamma^+ = \{(t,x); x > g(t)\}\setminus\gamma; \quad \gamma^- = \{(t,x); x < g(t)\}\setminus\gamma;$$
$$G = \gamma^+ \cap \theta^- \text{ or } G = \gamma^- \cap \theta^+ \ (G = \gamma^\pm \cap \theta^\mp);$$
$$\partial G = (\pm\gamma) \cup (\mp\theta).$$

For $k > 0$, denote $G_k^+ = (\gamma_k(a,b,\gamma))^\pm \cap (\gamma_k(a,b,\theta))^\mp$. Based on Theorem 2.3 of Guseynov (2016) for squares,
$$\int_{\partial G_k^+} f(t,x)dx = \int_{G_k^+} \frac{\partial}{\partial t} f(t,x)dtdx.$$

Since $\partial f/\partial t$ is summable and due to $g \in \mathcal{C}$, $|G_k^+ \setminus G| + |G \setminus G_k^+| \to 0$, $k \to \infty$, and we obtain formula (4.6). ∎

Proof of Theorem 2.5. With Lemma 4.3,
$$\int_{\gamma([0,s])} f(t,x)dx = \sum_k \sigma_k \int_{G_k} \frac{\partial}{\partial t} f(t,x)dtdx + \int_{\theta(s)} f(t,x)dx, 0 < s < 1,$$

(see notations in the statement of Theorem 2.5.) By Definition 2.1, the integral over the segment $\theta(s)$ is the line Riemann integral, which can be calculated as the definite Riemann integral
$$\int_{\theta(s)} f(t,x)dx = \frac{g(s)}{s} \int_0^s f\left(t, \frac{g(s)}{s}t\right)dt,$$

and we obtain the proof of (2.9).

To prove the integration by parts (2.11), we evaluate
$$\sum_k \sigma_k \int_{G_k} \frac{d}{dt} f(t)dtdx = -\sum_k \int_{a_k}^{b_k} \frac{d}{dt} f(t)dt \int_{\frac{g(s)}{s}t}^{g(t)} dx$$
$$= -\int_0^s \frac{d}{dt} f(t)dt \int_{\frac{g(s)}{s}t}^{g(t)} dx = -\int_0^s g(t) \frac{d}{dt} f(t)dt + \frac{g(s)}{s} \int_0^s \frac{d}{dt} f(t)tdt$$
$$= -\int_0^s g(t)df(t) + f(s)g(s) - \frac{g(s)}{s} \int_0^s f(t)dt,$$

since
$$\int_{[0,s]\setminus \cup_k(a_k,b_k)} \frac{d}{dt} f(t)dt \int_{\frac{g(s)}{s}t}^{g(t)} dx = 0.$$

Substituting this into (2.9) with $f(t,x) = f(t)$ yields (2.11). ∎

4. Two examples. If $g \in H_\alpha$ and $\beta + \alpha > 1, \beta, \alpha > 0$, then for $0 \le a < b \le 1, k_0 = [\log_2 2/|b-a|]$,
$$\mu(g,a,b,\beta) = \sum_{k > k_0} 2^{-\beta k} \sum_{[t_{k,n}, t_{k,n+1}] \subset [a,b]} |h_{k+1,2n}(g) - h_{k+1,2n+1}(g)|$$
$$\le |g, H_\alpha([0,T])| \sum_{k > k_0} 2^{-\beta k} 2^{-k\alpha} 2^{k-k_0} \le |g, H_\alpha| \frac{1}{1 - 2^{-(\alpha+\beta-1)}} |b-a|^{\alpha+\beta},$$

thus,
$$|g, G_{\beta,\alpha}([0,T])| = \sup_{0 \le a < b \le T} \frac{1}{(b-a)^{\beta+\alpha}} \mu(g,a,b,\beta) \le C|g, H_\alpha([0,T])|,$$

where $C$ depends only on $\beta + \alpha > 1$. Below, we show that, in general, $|x, G_{\beta,\gamma}([0,T])|$ and $|x, H_\alpha([0,T])|$, $x \in H_\alpha([0,T]) \cap G_{\beta,\gamma}([0,T])$, where $G_{\beta,\gamma}([0,T]) = \{g \in \mathcal{C}; |g, G_{\beta,\alpha}([0,T])| < \infty\}$, are not comparable if $\beta + \alpha < 1$, thus, both of them are independent characteristics of continuous functions.

**Example** 4.1. (High frequency oscillation input). For given $\beta, \alpha > 0, \beta + \alpha < 1, A \ge 0$, we construct a function $g \in H_\alpha$ such that $|g, H_\alpha| = 1$ and
$$C_1 2^{A(1-\alpha-\beta)/(1-\alpha)} \le |g, G_{\beta,\alpha\beta/(1-\alpha)}| \le C_2 2^{A(1-\alpha-\beta)/(1-\alpha)}, \tag{4.7}$$

where $C_1, C_2$ do not depend on $A$.

For $m > 1$ let $n_m$ be such that
$$\frac{A + \alpha m}{1 - \alpha} + 1 > n_m \geq \frac{A + \alpha m}{1 - \alpha}, \tag{4.8}$$
then
$$A + (1 - \alpha) \geq n_m - \alpha(n_m + m) \geq A.$$
Define the function $g_m: R \to R$,
$$g_m(t) = \sum_{p=0}^{2^{n_m}-1} f_{n_m+m,p}(t - 2^{-m}), \operatorname{supp}(g_m) \subset [2^{-m}, 2^{-m+1}],$$
where $f_{k,n}$ is the function defined in Theorem 4.1 with $\beta = \alpha$, and
$$g(t) = \sum_{m=1}^{\infty} g_m(t), t \in (0,1].$$
From (4.8) $n_m - \alpha(n_m + m) \geq A$,
$$\sum_{m=1}^{\infty} 2^{n_m} 2^{-\alpha(n_m+m)} = \infty,$$
and the function $g$ is unbounded variation on $[0,1]$.
If $k \leq n_m + m + 1$ then
$$A_k(g_m, a, b) = \sum_{[t_{k,n}, t_{k,n+1}] \subset [a,b]} |h_{k+1,2n}(g_m) - h_{k+1,2n+1}(g_m)| = 0, 0 \leq a < b \leq 1.$$
For $[a,b] \subset [2^{-m}, 2^{-m+1}], k_0 \geq m, k > n_m + m + 1$, since
$$2^{k-1-k_0} \leq n_k(b) - n_k(a) \leq 2^{k+1-k_0},$$
where $n_k(a), n_k(b)$ are from Definition 2.1 and $k_0 = [\log_2 2/|b - a|]$, then
$$A_k(g_m, a, b) = 2^{k+1} \sum_{[t_{k,n}, t_{k,n+1}] \subset [a,b]} \left| \int_{t_{k+1,2n}}^{t_{k+1,2n+1}} g_m(\tau) d\tau - \int_{t_{k+1,2n+1}}^{t_{k+1,2(n+1)}} g_m(\tau) d\tau \right|$$
$$\sim 2^{k+1} 2^{k-k_0} (2^{(1-\alpha)(n_m+m+1)} 2^{-k-1}) 2^{-k-1} = 2^{-\alpha} 2^{-k_0} 2^{(1-\alpha)(n_m+m)}$$
$$= 2^{-\alpha-k_0+m+n_m-\alpha(n_m+m)} \sim 2^{-k_0+m+A},$$
$$A_k(g_m, 2^{-m}, 2^{-m+1}) \sim 2^A.$$
($A \sim B$ means $C_1 A \leq B \leq C_2 A$ where $C_1, C_2$ are some constants.)
Having this, we evaluate
$$\mu(g, a, b, \beta) = \sum_{k > n_m+m+1} 2^{-\beta k} A_k(g_m, a, b) \sim 2^{-k_0+m+A} 2^{-\beta(n_m+m)} \sim 2^{-k_0+m+A-\beta(A+m)/(1-\alpha)}$$
$$= 2^{A(1-\alpha-\beta)/(1-\alpha)-k_0+m(1-\alpha-\beta)/(1-\alpha)} \leq 2^{\frac{A(1-\alpha-\beta)}{(1-\alpha)}} |b - a|^{\frac{\beta}{(1-\alpha)}},$$
and
$$\mu(g, 2^{-m}, 2^{-m+1}, \beta) \sim 2^{\frac{A(1-\alpha-\beta)}{(1-\alpha)}} 2^{-\frac{m\beta}{(1-\alpha)}}.$$
Another case that we should consider is $[a, b] \supset [2^{-i}, 2^{-i+1}]$, where $[2^{-i}, 2^{-i+1}]$ is the biggest segment that $[a, b]$ contains. For the lower bound, based on previous estimate, we obtain
$$\mu(g, a, b, \beta) \geq \mu(g, 2^{-i}, 2^{-i+1}, \beta) \sim 2^{\frac{A(1-\alpha-\beta)}{(1-\alpha)}} 2^{-\frac{i\beta}{(1-\alpha)}} \geq C 2^{\frac{A(1-\alpha-\beta)}{(1-\alpha)}} |b - a|^{\frac{\beta}{(1-\alpha)}},$$
where $C$ does not depend on $a, b, A$. For the upper bound, consider
$$\mu(g, a, b, \beta) \leq \mu(g, 0, 2^{-i+2}, \beta) = \sum_{k > i-2} 2^{-\beta k} \sum_{m, n_m+m+1 < k} A_k(g_m).$$
If $(A + m)/(1 - \alpha) + 1 < k$ then $m \leq (1 - \alpha)k - A - (1 - \alpha)$ and $(A + m)/(1 - \alpha) + 1 < n_m + m + 1$. Thus,
$$\mu(g, a, b, \beta) \leq \sum_{k > i-2} 2^{-\beta k} \sum_{m=i-1}^{(1-\alpha)k-A} 2^A \leq 2^A \sum_{k \geq \frac{i+A-1}{1-\alpha}} 2^{-\beta k}((1-\alpha)k - A - i + 2)$$

$$= 2^A \sum_{m=1}^{\infty} m \, 2^{-\beta(m+A+i-2)/(1-\alpha)} \leq C 2^{\frac{A(1-\alpha-\beta)}{(1-\alpha)}} |b-a|^{\frac{\beta}{(1-\alpha)}},$$

where $C$ does not depend on $a, b, A$. ∎

Let $g^A$ be the function constructed in Example 4.1 for $A > 0, |g^A, H_\alpha| = 1$, then
$$\infty > |g^A, G_{\beta,\alpha\beta/(1-\alpha)}| \geq C 2^{A(1-\alpha-\beta)/(1-\alpha)} |g^A, H_\alpha|,$$
where $C$ does not depend on A. This shows that for $\beta + \alpha < 1$, $|g^A, G_{\beta,\alpha\beta/(1-\alpha)}| \to \infty$ when $A \to \infty$, but $|g^A, H_\alpha|$ remains constant which means $|g^A, G_{\beta,\alpha\beta/(1-\alpha)}|$ and $|g^A, H_\alpha|$ are not comparable.

For $A \geq 0$ we obtain for any $\beta, \alpha > 0, \beta + \alpha < 1$, that $|g^A, G_{\beta,\alpha\beta/(1-\alpha)}| < \infty$ which proves that the set $G_{\beta,\alpha\beta/(1-\alpha)} \cap H_\alpha$ is not empty.

By Wiener, for any $\alpha \in (0, 2^{-1}), W(H_\alpha) = 1$, and by Young's theorem, for $\beta > 2^{-1}, H_\alpha \subset GB_{\beta,\alpha}$ for $\alpha = 2^{-1} - 2^{-1}(\beta - 2^{-1})$, thus, $W(GB_{\beta,\alpha}) = 1$ or the set
$$GB_{\beta,\alpha} = \{g \in H_\alpha; \forall f \in CH_\beta(R_g) \, \forall a \, \forall b, 0 \leq a < b \leq 1, \exists \int_a^b f(t, g(t)) dg(t)\}$$
has full Wiener measure. For $\alpha + \beta < 1, \beta, \alpha > 0$, in the example below, we construct a function $g \in H_\alpha$ such that integral (2.1) over $g$ cannot be applied to the whole class $H_\beta$ ($g \notin GB_{\beta,\alpha}$). This indicates that the Wiener measure of class $H_\alpha$ cannot be used to evaluate $GB_{\beta,\alpha}$ for $\beta \leq 2^{-1}$, similar to the case of $\beta > 2^{-1}$. Meanwhile, based on Theorem 2.1, there is a subclass $H_\alpha \cap G_{\beta,\gamma}$ which is not empty, see previous example, such that integral (2.1) over $g \in H_\alpha \cap G_{\beta,\gamma}$ exists for all functions from $H_\beta$.

**Example 4.2** ($g \in H_\alpha \setminus GB_{\beta,\alpha}$). Let $\gamma = \alpha + \beta < 1, \beta, \alpha > 0, a_n = (2^{1-\gamma} - 1)2^{-n(1-\gamma)}$,
$$\sum_{n=1}^{\infty} a_n = 1, J_k = \left[\sum_{n=k+1}^{\infty} a_n, \sum_{n=k}^{\infty} a_n\right], k > k_0 = \left[\frac{1}{\gamma} \log_2 \frac{4}{2^{1-\gamma}-1} + 1\right].$$
For each segment $I_{k,m} = [2^{-k}m, 2^{-k}(m+1)] \subset J_k, m = 0, \ldots, 2^k - 1$, we replace line segment $[2^{-k}m, 2^{-(k+1)}(2m+1)]$ with two sides of the triangle with the vertex $(2^{-(k+2)}(4m+1), 2^{-(k+1)\alpha})$. For the number $N_k$ of segments $I_{k,m} \subset J_k, (2^{1-\gamma} - 1)2^{k\gamma} \geq N_k \geq (2^{1-\gamma} - 1)2^{k\gamma} - 3$, and for $k > k_0, N_k > 1$. The resulting Jordan curve can be parameterized by a function $g(t), t \in [0,1]$, from $H_\alpha$, and since
$$\sum_{k=k_0}^{\infty} 2^{-(k+1)\alpha} N_k = \infty,$$
it is unbounded variation. Additionally, the estimate of $|B_{k+1}(g)|$ on the segment $J_k$ is
$$|B_{k+1}(g)| \geq N_k 2^{-(k+1)\alpha} 2^{-(k+2)} \geq 2^{-(2+\alpha)}(2^{1-\gamma} - 1) 2^{k(\beta-1)} - 3 * 2^{-(k+2)(\alpha+1)},$$
and condition (2.2)
$$\sum_{k=k_0}^{k_1} 2^{k(1-\beta)} |B_{k-1}(g)|$$
$$\geq 2^{-2-\alpha}(2^{1-\gamma} - 1) \sum_{k=k_0}^{k_1} 2^{k(1-\beta)} 2^{-(k-2)(1-\beta)} - 3 \sum_{k=k_0}^{k_1} 2^{-(k+1)(\alpha+1)-1} \to \infty, k_1 \to \infty.$$
Thus, based on Theorem 2.2 (necessity), integral (2.1) with respect to the function $g \in H_\alpha$ cannot be applied to the whole class $H_\beta$ if $\gamma = \alpha + \beta < 1$. ∎